\documentclass[amstex,12pt, amssymb]{article}

\usepackage{mathtext}
\usepackage[cp1251]{inputenc}
\usepackage[T2A]{fontenc}
\usepackage[dvips]{graphicx}
\usepackage{amsmath}
\usepackage{amssymb}
\usepackage{amsxtra}
\usepackage{latexsym}
\usepackage{ifthen}

\def\Xint#1{\mathchoice
   {\XXint\displaystyle\textstyle{#1}}
   {\XXint\textstyle\scriptstyle{#1}}
   {\XXint\scriptstyle\scriptscriptstyle{#1}}
   {\XXint\scriptscriptstyle\scriptscriptstyle{#1}}
   \!\int}
\def\XXint#1#2#3{{\setbox0=\hbox{$#1{#2#3}{\int}$}
     \vcenter{\hbox{$#2#3$}}\kern-.5\wd0}}
\def\dashint{\Xint-}

\newcounter{lemma}[section]

\newcounter{corollary}[section]

\newcounter{remark}[section]

\newcounter{theorem}[section]

\newcounter{proposition}[section]

\newcounter{example}

\numberwithin{equation}{section}

\begin{document}

\markboth{\centerline{V.~GUTLYANSKII, V.~RYAZANOV, E.~SEVOST'YANOV
and E.~YAKUBOV}}{\centerline{BELTRAMI EQUATIONS AND ASYMPTOTIC
HOMOGENEITY}}

\def\cc{\setcounter{equation}{0}
\setcounter{figure}{0}\setcounter{table}{0}}

\overfullrule=0pt


\author{Gutlyanskii V., Ryazanov V., Sevos'yanov E., Yakubov E.}

\title{Beltrami equations and mappings\\
with asymptotic homogeneity at infinity}

\date{\today}
\maketitle

\begin{abstract}
First of all, we prove that the BMO condition by John-Nirenberg
leads in the natural way to the asymptotic homogeneity at the origin
of regular homeomorphic solutions of the degenerate Beltrami
equations. Then on this basis we establish a series of criteria for
the existence of regular homeomorphic solutions of the degenerate
Beltrami equations in the whole complex plane with asymptotic
homogeneity at infinity. These results can be applied to the fluid
mechanics in strictly anisotropic and inhomogeneous media because
the Beltrami equation is a complex form of the main equation of
hydromechanics.
\end{abstract}





\bigskip

{\bf Dedicated to the memory of mathematicians Fritz John and Louis
Nirenberg}


\section{Introduction}

A real-valued function $u$ in a domain $D$ in ${\Bbb C}$ is said to
be of {\bf bounded mean oscillation} in $D$, abbr. $u\in{\rm
BMO}(D)$, if $u\in L_{\rm loc}^1(D)$ and
\begin{equation}\label{lasibm_2.2_1}\Vert u\Vert_{*}:=
\sup\limits_{B}{\frac{1}{|B|}}\int\limits_{B}|u(z)-u_{B}|\,dm(z)<\infty\,,\end{equation}
where the supremum is taken over all discs $B$ in $D$ and
$$u_{B}={\frac{1}{|B|}}\int\limits_{B}u(z)\,dm(z)\,.$$

Recall that the class BMO was introduced by John and Nirenberg
(1961) in the paper \cite{JN} and soon became an important concept
in harmonic analysis, partial differential equations and related
areas, see e.g. \cite{HKM} and \cite{RR}.


A function $\varphi$ in BMO is said to have {\bf vanishing mean
oscillation}, abbr. $\varphi\in{\rm VMO}$, if the supremum in
(\ref{lasibm_2.2_1}) taken over all balls $B$ in $D$ with
$|B|<\varepsilon$ converges to $0$ as $\varepsilon\to0$. Recall that
VMO has been introduced by Sarason in \cite{Sar}. There are a number
of papers devoted to the study of partial differential equations
with coefficients of the class VMO, see e.g. \cite{CFL}, \cite{IS},
\cite{MRV}, \cite{Pal}, \cite{Ra$_1$} and \cite{Ra$_2$}. Note by the
way that $W^{\,1,2}\left({{D}}\right) \subset VMO
\left({{D}}\right),$ see \cite{BN}.

Let $D$ be a domain in the complex plane ${\Bbb C}$, i.e., a
connected open subset of ${\Bbb C}$, and let $\mu:D\to{\Bbb C}$ be a
measurable function with $|\mu(z)|<1$ a.e. (almost everywhere) in
$D$. A {\bf Beltrami equation} is an equation of the form
\begin{equation}\label{eqBeltrami} \overline{\partial}f(z)=\mu(z)\cdot\partial f(z)\,\end{equation}
with the formal complex derivatives
$\overline{\partial}f=(f_x+if_y)/2$, $\partial f=(f_x-if_y)/2$,
$z=x+iy$, where $f_x$ and $f_y$ are partial derivatives of $f$ in
$x$ and $y$, correspondingly. The function $\mu$ is  said to be the
{\bf complex coefficient} and
\begin{equation}\label{eqKPRS1.1}K_{\mu}(z)\ :=\ \frac{1+|\mu(z)|}{1-|\mu(z)|}\end{equation}
the {\bf dilatation quotient} of the equation (\ref{eqBeltrami}).
The Beltrami equation is called {\bf degenerate} if ${\rm ess}\,{\rm
sup}\,K_{\mu}(z)=\infty$. Homeomorphic solutions of the Beltrami
equations with $K_{\mu}\le Q<\infty$ in the Sobolev class
$W^{1,1}_{\rm loc}$ are called {\bf Q-quasiconformal mappings}.

It is well known that if $K_{\mu}$ is bounded, then the Beltrami
equation has homeomorphic solutions, see e.g. \cite{Ahl},
\cite{Bojar}, \cite{LV} and \cite{Vek}. Recently, a series of
effective criteria for the existence of homeomorphic $W^{1,1}_{\rm
loc}$ solutions have been also established for degenerate Beltrami
equations, see e.g. historic comments with relevant references in
monographs the \cite{AIM}, \cite{GRSY} and \cite{MRSY}, in
BMO-article \cite{RSY} and in the surveys \cite{GRSY$_*$} and
\cite{SY}.

These criteria were formulated both in terms of $K_{\mu}$ and the
more refined quantity that takes into account not only the modulus
of the complex coefficient $\mu$ but also its argument
\begin{equation}\label{eqTangent} K^T_{\mu}(z,z_0)\ :=\
\frac{\left|1-\frac{\overline{z-z_0}}{z-z_0}\mu (z)\right|^2}{1-|\mu
(z)|^2} \end{equation} that is called the {\bf tangent dilatation
quotient} of the Beltrami equation with respect to a point
$z_0\in\mathbb C$, see e.g. \cite{And}, \cite{BGR$_1$},
\cite{BGR$_2$}, \cite{GMSV}, \cite{Le} and
\cite{RSY}-\cite{RSY$_5$}. Note that
\begin{equation}\label{eqConnect} K^{-1}_{\mu}(z)\leqslant K^T_{\mu}(z,z_0) \leqslant K_{\mu}(z)
\ \ \ \ \ \ \ \forall\ z\in D\, ,\ z_0\in \Bbb C\ .\end{equation}
The geometrical sense of $K^T_{\mu}$ can be found e.g. in the
monographs \cite{GRSY} and \cite{MRSY}.

\medskip

A function $f$ in the Sobolev class $W^{1,1}_{\rm loc}$ is called a
{\bf regular solution} of the Beltrami equation (\ref{eqBeltrami})
if $f$ satisfies it a.e.  and its Jacobian
$J_f(z)=|{\partial}f(z)|^2-|\overline{\partial}f(z)|^2>0$ a.e. in
$\mathbb C.$
 The notion of such a solution was probably first
introduced in \cite{BGR}.

\bigskip

By the well-known Gehring-Lehto-Menchoff theorem, see \cite{GeLe}
and \cite{Me}, see also the monographs \cite{Ahl} and \cite{LV},
each homeomorphic $W^{1,1}_{\rm loc}$ solution $f$ of the Beltrami
equation is differentiable a.e. Recall that a function $f:
D\to\mathbb C$ is {\bf differentiable by Darboux–Stolz  at a point}
$z_0\in D$ if
\begin{equation}\label{eqDaSt}
f(z) - f(z_0) = {\partial}f(z_0)\cdot (z-z_0) +
\overline{\partial}f(z_0)\cdot \overline{(z-z_0)} + o(|z-z_0|)
\end{equation} where $o(|z-z_0|)/|z-z_0|\to 0$ as $z\to z_0$.
Moreover, $f$ is called {\bf conformal at the point} $z_0$ if in
addition $f_{\overline{z}} (z_0) =0$ but $f_z (z_0) \neq 0$.

The example $w = z(1-\ln |z|)$ of B.V. Shabat, see \cite{Be}, p. 40,
shows that, for a continuous complex characteristic $\mu(z)$, the
quasiconformal mapping $w = f(z)$ can be non-differentiable  by
Darboux–Stolz  at the origin. If the characteristic $\mu(z)$ is
continuous at a point $z_0\in D$, then, as was first established,
apparently, by P.P. Belinskij in \cite{Be}, p. 41, the mapping $w =
f(z)$ is differentiable at $z_0$ in the following meaning:
\begin{equation}\label{eqDifBel}
\Delta w = A(\rho) \left[ \Delta z + \mu_0 \overline{\Delta z} + o
(\rho) \right] ,
\end{equation}
where $\mu_0 = \mu(z_0), \ \rho = |\Delta z + \mu_0 \overline{\Delta
z}|, \ A(\rho)$ depends only on $\rho$ and $o(\rho) / \rho \to 0$ as
$\rho \to 0$. As it was clarified later in \cite{R}, see also
\cite{GR}, here $A(\rho)$ may not have a limit with $\rho\to 0$,
however,
\begin{equation}\label{eqDifBelR}
\lim\limits_{\rho \to 0} \frac{A(t \rho)}{A(\rho)} = 1 \quad
    \forall\,t>0\ .
\end{equation}

Following \cite{R}, a mapping $f:D\to\mathbb C$ is called {\bf
differentiable by Belinskij at a point} $z_0\in D$ if conditions
(\ref{eqDifBel}) and (\ref{eqDifBelR}) hold with some
$\mu_0\in\mathbb D:= \{ \mu\in\mathbb C: |\mu|<1 \} $. Note that
here, in the case of discontinuous $\mu(z)$, it is not necessary
$\mu_0=\mu(z_0)$. If in addition $\mu_0=0$, then $f$ is called {\bf
conformal by Belinskij at the point} $z_0$.

For quasiconformal mappings $f:D\to\mathbb C$ with $f(0)=0\in D$, it
was shown in \cite{R}, see also \cite{GR}, that the conformality by
Belinskij of $f$ at the origin is equivalent to each of its
properties:
\begin{equation}\label{eq_4.5}
    \boxed{\ \lim\limits_{\tau \to 0} \frac{f(\tau \zeta)}{f(\tau)} = \zeta\ }\
    \ \ \ \ \ \ \mbox{along the ray $\tau > 0$}\ \ \ \ \ \forall\ \zeta \in
    \mathbb{C}\ ,
\end{equation}
\begin{equation}\label{eq_4.6}
 \boxed{\ \lim\limits_{z \to 0} \left\{ \frac{f(z')}{f(z)} - \frac{z'}{z} \right\}=
 0\ }\ \ \ \ \ \ \mbox{along $z, z' \in
\mathbb{C}$, $|z'| < \delta |z|$,}\ \ \ \ \ \forall\ \delta>0\ ,
\end{equation}
\begin{equation}\label{eq_4.7}
  \boxed{\ \lim\limits_{z \to 0} \frac{f(z \zeta)}{f(z)} = \zeta\ }
  \ \ \ \ \ \ \ \ \ \mbox{along $z \in \mathbb{C}^{*}:=
\mathbb{C} \setminus \{0\}$}\ \ \ \ \ \forall\ \zeta \in \mathbb{C}\
,
\end{equation}
and, finally, to the property of the limit in (\ref{eq_4.7}) to be
locally uniform with respect to $\zeta\in\mathbb C$.

Following the article \cite{GR}, the property (\ref{eq_4.7}) of a
mapping $f:D\to\mathbb C$ with $f(0)=0\in D$ is called its {\bf
asymptotic homogeneity} at $0$. In the sequel we sometimes write
(\ref{eq_4.7}) in the shorter form $f(\zeta z)\sim\zeta f(z)$.

In particular, we obtain from (\ref{eq_4.6}) under $|z^{\prime}| =
|z|$ that
\begin{equation}\label{eq_4.8}
    \boxed{\ \lim\limits_{r \to 0}\ \frac{\max\limits_{|z| =r} |f(z)|} {\min\limits_{|z| =r} |f(z)|}\ =\
    1\ }
\end{equation}
i.e., that the Lavrent'iev characteristic is equal $1$ at the
origin. It is natural to say in the case of (\ref{eq_4.8}) that the
mapping $f$ is {\bf conformal by Lavrent'iev} at $0$. As we see, the
usual conformality implies the conformality by Belinskij and the
latter implies the conformality by  Lavrent'iev at the origin
meaning geometrically that the infinitesimal circle centered at zero
is transformed into an infinitesimal circle also centered at zero.

However, condition (\ref{eq_4.7}) much more stronger than condition
(\ref{eq_4.8}). We obtain from (\ref{eq_4.7}) also asymptotic
preserving angles
\begin{equation}\label{eq_4.9}
    \boxed{\ \lim\limits_{z \to 0}\ \arg\, \left[\frac{f(z \zeta)}{f(z)} \right] = \arg
    \zeta\ }\ \ \ \ \ \forall\ \zeta \in \mathbb{C}^*
\end{equation}
and asymptotic preserving moduli of infinitesimal rings
\begin{equation}\label{eq_4.10}
    \boxed{\ \lim\limits_{z \to 0}\ \frac{|f(z \,\zeta)|}{|f(z)|}\ =\
    |\zeta|\ }\ \ \ \ \ \forall\ \zeta \in \mathbb{C}^*\ .
\end{equation}
The latter two geometric properties characterize asymptotic
homogeneity and demonstrate that it is close to the usual
conformality.

It should be noted that, despite (\ref{eq_4.10}), an asymptotically
homogeneous map can send radial lines to infinitely winding spirals,
as shown by the example $f(z) = z e^{i \sqrt{- \ln |z|}}$, see
\cite{Be}, p. 41. Moreover, the above Shabat example shows that the
conformality by Belinskij admits infinitely great tensions and
pressures at the corresponding points.

It was shown in \cite{R} that a quasiconformal mapping
$f:D\to\mathbb C$, whose complex cha\-rac\-te\-ris\-tic $\mu(z)$ is
approximately continuous at a point $z_0\in D$, is differentiable by
Belinskij at the point with $\mu_0=\mu(z_0)$ and, in particular, is
asymptotically homogeneous if $\mu(z_0)=0$. Recall that $\mu(z)$ is
called {\bf approximately continuous at the point} $z_0$ if there is
a measurable set $E$ such that $\mu(z)\to\mu(z_0)$ as $z\to z_0$ in
$E$ and $z_0$ is a point of density for $E$, i.e.,
$$
\lim\limits_{\varepsilon \to 0}\ \frac{|E \cap D(z_0, \varepsilon)|}
{|D (z_0, \varepsilon)|}\ =\ 1\ ,
$$
where $D(z_0, \varepsilon) = \{ z \in \mathbb{C} : |z-z_0| \leq
\varepsilon \}$. Note also that, for functions $\mu$ in
$L^{\infty}$, the points of approximate continuity coincide with the
Lebesgue points of $\mu$, i.e., such $z_0$ for which
$$
\lim\limits_{r\to 0}\ \frac{1}{r^2}\int\limits_{|z-z_0|<r}\,
|\mu(z)\ -\ \mu(z_0)|\ dm(z)\ =\ 0\ ,
$$
where $dm(z):=dxdy$, $z=x+iy$, stands to the Lebesgue measure (area)
in $\mathbb C$.

In comparison with the previous version in arXiv, the above results
on the asymptotic homogeneity, i.e., on the conformality by
Belinskij, are extended to the degenetare Beltrami equations with
majorants of its dilatation $K_{\mu}$ in BMO.

\section{FMO and the main lemma with participation of BMO}

Here and later on, we apply the notations
$$
\mathbb D(z_0,r)\ :=\ \{ z\in\mathbb C: |z-z_0|\ <\ r\}\ ,\ \ \ \ \
\mathbb D(r)\ :=\ \mathbb D(0,r) \ ,\ \ \ \ \ \mathbb D\ :=\ \mathbb
D(0,1)\ ,
$$
and of the mean value of integrable functions $\varphi$ over the
disks $\mathbb D(z_0,r)$
$$\dashint_{\mathbb D(z_0,r)}\varphi(z)\,
dm(z)\ :=\ \frac{1}{|\mathbb D(z_0,r)|}\int\limits_{\mathbb
D(z_0,r)}\varphi(z)\, dm(z)\ .$$

Following \cite{IR}, we say that a function $\varphi:D\to{\Bbb R}$
has {\bf finite mean oscillation} at a point $z_0\in D$, abbr.
$\varphi\in{\rm FMO}(z_0)$, if
\begin{equation}\label{FMO_eq2.4}\overline{\lim\limits_{\varepsilon\to0}}\ \ \
\dashint_{\mathbb
D(z_0,\varepsilon)}|{\varphi}(z)-\widetilde{\varphi}_{\varepsilon}(z_0)|\,dm(z)<\infty\,,\end{equation}
where \begin{equation}\label{FMO_eq2.5}
\widetilde{\varphi}_{\varepsilon}(z_0)=\dashint_{\mathbb
D(z_0,\varepsilon)} {\varphi}(z)\,dm(z)\ .\end{equation} Note that
the condition (\ref{FMO_eq2.4}) includes the assumption that
$\varphi$ is integrable in some neighborhood of the point $z_0$. We
say also that a function $\varphi:D\to{\Bbb R}$ is of {\bf finite
mean oscillation in D}, abbr. $\varphi\in{\rm FMO(D)}$ or simply
$\varphi\in{\rm FMO}$, if $\varphi\in{\rm FMO}(z_0)$ for all points
$z_0\in D$.

\medskip

{\bf Remark 1.} It is evident that ${\rm BMO(D)}\subset{\rm
BMO(D)}_{\rm loc}\subset{\rm FMO(D)}$ and it is well-known by the
John-Nirenberg lemma that ${\rm BMO}_{\rm loc}\subset L_{\rm loc}^p$
for all $p\in[1,\infty)$, see e.g. \cite{JN} or \cite{RR}. However,
FMO is not a subclass of $L_{\rm loc}^p$ for any $p>1$ but only of
$L_{\rm loc}^1$, see e.g. example 2.3.1 in \cite{GRSY} or the
example in \cite{RSY$_6$}. Thus, the class FMO is much more wider
than ${\rm BMO}_{\rm loc}$.

\medskip

The following statement is obvious by the triangle inequality.

\medskip

{\bf Proposition 1.} {\it If, for a  collection of numbers
$\varphi_{\varepsilon}\in{\Bbb R}$,
$\varepsilon\in(0,\varepsilon_0]$,
\begin{equation}\label{FMO_eq2.7}\overline{\lim\limits_{\varepsilon\to0}}\ \ \
\dashint_{\mathbb
D(z_0,\varepsilon)}|\varphi(z)-\varphi_{\varepsilon}|\,dm(z)<\infty\,,\end{equation}
then $\varphi $ is of finite mean oscillation at $z_0$.}

\medskip

In particular, choosing here  $\varphi_{\varepsilon}\equiv0$,
$\varepsilon\in(0,\varepsilon_0]$ in Proposition 1, we obtain the
following.

\medskip

{\bf Corollary 1.} {\it If, for a point $z_0\in D$,
\begin{equation}\label{FMO_eq2.8}\overline{\lim\limits_{\varepsilon\to 0}}\ \ \
\dashint_{\mathbb D(z_0,\varepsilon)}|\varphi(z)|\,dm(z)<\infty\,,
\end{equation} then $\varphi$ has finite mean oscillation at
$z_0$.}

\medskip

Recall that a point $z_0\in D$ is called a {\bf Lebesgue point} of a
function $\varphi:D\to{\Bbb R}$ if $\varphi$ is integrable in a
neighborhood of $z_0$ and \begin{equation}\label{FMO_eq2.7a}
\lim\limits_{\varepsilon\to 0}\ \ \ \dashint_{\mathbb
D(z_0,\varepsilon)}
|\varphi(z)-\varphi(z_0)|\,dm(z)=0\,.\end{equation} It is known
that, almost every point in $D$ is a Lebesgue point for every
function $\varphi\in L^1(D)$. Thus, we have by Proposition 1 the
next corollary.

\medskip

{\bf Corollary 2.} {\it Every locally integrable function
$\varphi:D\to{\Bbb R}$ has a finite mean oscillation at almost every
point in $D$.}

\medskip

{\bf Remark 2.} Note that the function
$\varphi(z)=\log\left(1/|z|\right)$ belongs to BMO in the unit disk
$\mathbb D$, see e.g. \cite{RR}, p. 5, and hence also to FMO.
However, $\widetilde{\varphi}_{\varepsilon}(0)\to\infty$ as
$\varepsilon\to0$, showing that condition (\ref{FMO_eq2.8}) is only
sufficient but not necessary for a function $\varphi$ to be of
finite mean oscillation at $z_0$.

\medskip

Versions of the next statement has been first proved for the class
BMO in \cite{RSY}. For the FMO case, see the papers \cite{IR, RS,
RSY$_2$, RSY$_3$}  and the monographs \cite{GRSY} and \cite{MRSY}.
Here we prefer to use its following version, see Lemma 2.1 in
\cite{RSY$_6$}, cf. also Lemma 5.3 in the monograph \cite{GRSY}:

\medskip

{\bf Proposition 2.} {\it  Let ${\varphi}: D\rightarrow
{\mathbb{R}}$ be a nonnegative function with finite mean oscillation
at $0\in {D}$ and integrable in the disk $\mathbb D(1/2)\subset D.$
Then
\begin{equation}  \label{FMO_eq2.9}
\int\limits_{A(\varepsilon, 1/2)}\frac{{\varphi} (z)\, dm(z)}
{\left(|z| \log_2 \frac{1}{|z|}\right)^2} \le\ C \cdot \log_2 \log_2
\frac{1}{\varepsilon}\ \ \ \ \ \ \ \forall\ {\varepsilon}\in
(0,1/4)\ ,
\end{equation}
where
\begin{equation}  \label{FMO_eq2.10a}
C = 4\pi\ ({\varphi} _0 + 6d_0)\;,
\end{equation}
${\varphi} _0$ is the average of ${\varphi}$ over the disk $\mathbb
D(1/2)$ and $d_0$ is the maximal dispersion of ${\varphi}$ in
$\mathbb D(1/2).$ }

\medskip

Recall that the \textbf{maximal dispersion} of the function
${\varphi}$ in the disk $\mathbb D(z_0, r_0)$ is the quantity
\begin{equation}  \label{FMO_eq2.4c}
\sup\limits_{r\in(0,r_0]}\ \ \ \mathchoice
{{\setbox0=\hbox{$\displaystyle{\textstyle -}{\int}$}
\vcenter{\hbox{$\textstyle -$}}\kern-.5\wd0}}
{{\setbox0=\hbox{$\textstyle{\scriptstyle -}{\int}$}
\vcenter{\hbox{$\scriptstyle -$}}\kern-.5\wd0}}
{{\setbox0=\hbox{$\scriptstyle{\scriptscriptstyle -}{\int}$}
\vcenter{\hbox{$\scriptscriptstyle -$}}\kern-.5\wd0}}
{{\setbox0=\hbox{$\scriptscriptstyle{\scriptscriptstyle -}{\int}$}
\vcenter{\hbox{$\scriptscriptstyle -$}}\kern-.5\wd0}}
\!\int_{\mathbb D( z_0 ,r)}
|{\varphi}(z)-\widetilde{{\varphi}}_r(z_0)|\ dm(z)\ .
\end{equation}
Here and later on, we also use the following designations for the
spherical rings in $\mathbb C\ $:
\begin{equation}  \label{FMO_eq2.10}
A(z_0,r_1,r_2)\ :=\ \{z\in{\mathbb{C}}:\ r_1\, <\, |z-z_0|\, <\,
r_2\} \ ,\ \ \ \ A(r_1,r_2)\ :=\ A(0,r_1,r_2)\ .
\end{equation}

\medskip

Further, we denote by $M$ the conformal modulus (or $2-$modulus) of
a family of paths in $\mathbb C$, see e.g. \cite{Va}. Moreover,
given sets $E$ and $F$ and a domain $D$ in ${{\Bbb C}}$, we denote
by $\Gamma(E, F, D)$ the family of all paths $\gamma:[0,
1]\rightarrow{{\Bbb C}}$ joining $E$ and $F$ in $D,$ that is,
$\gamma(0)\in E,$ $\gamma(1)\in F$ and $\gamma(t)\in D$ for all
$t\in (0, 1).$

\medskip

Let $Q:{\Bbb C}\rightarrow (0,\infty)$ be a Lebesgue measurable
function. A mapping $f:D\rightarrow{{\Bbb C}}$ is called a {\bf ring
Q}$-${\bf mapping at a point} $z_0\in D$, if
\begin{equation}\label{eq3*!gl0} M(f(\Gamma(\mathbb S(z_0, r_1),\,\mathbb S(z_0, r_2),
D)))\ \leqslant\ \int\limits_A Q(z)\cdot \eta^2(|z-z_0|)\, dm(z)
\end{equation}
for each spherical ring $A=A(z_0,r_1,r_2)$ with arbitrary
$0<r_1<r_2<\delta_0:={\rm dist}\,(z_0,
\partial D)$ and all Lebesgue measurable functions $\eta:(r_1,
r_2)\rightarrow [0, \infty]$ such that
\begin{equation}\label{eq*3gl0}
\int\limits_{r_1}^{r_2}\eta(r)\,dr\ \geqslant\ 1\ .
\end{equation}
Here we use also the notations for the circles in $\mathbb C$
centered at a point $z_0$
$$
\mathbb S(z_0,r_0)\ =\ \{\, z\in\mathbb C:\ |z\,-\,z_0|\,=\,r_0\,\}\
.
$$


{\bf Remark 3.} Recall that regular homeomorphic solutions of the
Beltrami equation (\ref{eqBeltrami}) are $Q_{z_0}-$map\-pings with
$Q_{z_0}(z)=K^T_{\mu}(z,z_0)$ and, in particular, $Q-$mappings with
$Q(z)=K_{\mu}(z)$ at each point $z_0\in D$, see \cite{Sal}, see also
Theorem 2.2 in \cite{GRSY}, cf. Theorem 3.1 in \cite{LSS} and
Theorem 3 in \cite{KPRS}.

\bigskip

Later on, in the extended complex plane $\overline{{{\Bbb
C}}}={{\Bbb C}}\cup\{\infty\},$ we use the {\bf spherical (chordal)
metric} $s$ defined by the equalities
\begin{equation}\label{eq3C}
s(z,\zeta)=\frac{|z-\zeta|}{\sqrt{1+{|z|}^2}
\sqrt{1+{|\zeta|}^2}}\,,\quad z\ne \infty\ne \zeta\,,
\quad\,s(z,\infty)=\frac{1}{\sqrt{1+{|z|}^2}}\ ,
\end{equation}
see e.g.~\cite[Definition~12.1]{Va}. For a given set $E$ in
$\overline{{\Bbb C}},$ we also use its {\bf spherical diameter}
\begin{equation}\label{eq47***}
s(E):=\sup\limits_{z,\zeta\in E}s(z, \zeta)\,.
\end{equation}

Given a domain $D$ in $\mathbb C$, a prescribed point $z_0\in D$ and
a measurable $Q:D\to(0,\infty),$ later on $\frak{R}^{\Delta}_Q$
denotes the class of all ring $Q-$homeomorphisms $f$ at $z_0$ in $D$
with
$$s(\overline{{\mathbb C}}\setminus f(D))\ \geq\ \Delta>0\ .$$

\medskip

The following statement, see Theorem 4.3 in \cite{RSY$_6$}, provides
us by the effective estimates of the distortion of the spherical
distance under the ring $Q-$homeomorphisms, and it follows just on
the basis of Proposition 2 on FMO functions above.

\medskip

{\bf Proposition 3.} {\it  Let $f\in\frak{R}^{\Delta}_Q(D)$ with
$\Delta>0$ and $Q: D\rightarrow {\mathbb{R}}$ be a nonnegative
function with finite mean oscillation at $\zeta_0\in {D}$ and
integrable in the disk $\mathbb D(\zeta_0,\varepsilon_0)\subset D$,
$\varepsilon_0>0$. Then
\begin{equation}  \label{FMO_eq3.1m}
s(f(\zeta), f(\zeta_0))\ \le\ \frac{32}{\Delta}\cdot {\left(\log \frac{2\varepsilon_0}{|\zeta-\zeta_0|}%
\right)^{-\frac{1}{\alpha_0}}}\ \ \ \ \ \ \ \ \ \ \forall\ \zeta\in
\mathbb D(\zeta_0,\varepsilon_0/2)\ ,
\end{equation}
where
\begin{equation}  \label{FMO_eq2.10a}
\alpha_0\ =\ 2(q_0 + 6d_0)\ ,
\end{equation}
$q_0$ is the average of $Q$ over $\mathbb D(\zeta_0,\varepsilon_0)$
and $d_0$ is the maximal dispersion of $Q$ in $\mathbb
D(\zeta_0,\varepsilon_0).$ }

\medskip

Propositions 2 and 3 are key in establishing equicontinuity of
classes of mappings associated with asymptotic homogeneity in the
proof of the central lemma involving BMO.

\medskip

{\bf Lemma 1.} {\it Let $D$ be a domain in $\mathbb{C}$, $0\in D$,
and let $f:D\rightarrow{\mathbb{C}}$ be a regular homeomorphic
solution of the Beltrami equation (\ref{eqBeltrami}) with $f(0)=0$.
Suppose that its dilatation $K_{\mu}$ has a majorant $Q\in$ {\rm
BMO(D)}. Then the family of mappings $f_z(\zeta) := f(\zeta \, z) /
f(z)$ is equicontinuous with respect to the spherical metric at each
point $\zeta_0\in\mathbb C$ as $z\to 0$ along $z \in
\mathbb{C}^{*}:= \mathbb{C} \setminus \{0\}$.}

\medskip

\begin{proof} Indeed, for $\zeta_0\in\mathbb D(\delta)$, $\delta>1$, $0<\delta_*<{\rm dist}\, (0,\partial
D)$, $\tau_*:=\delta_*/\delta<\delta_*$, we see that
$$
\mathbb D(z\zeta_0,\rho_{z})\ \subseteq\ \mathbb
D(\delta_*)\subseteq D\ , \ \ \   \mbox{where $\rho_{z}:=\delta_*
-|z\zeta_0|\geq \delta_*(1-|\zeta_0|/\delta)>0$,\
$z\in\overline{\mathbb D(\tau_*)}\setminus\{ 0\}$\, .}
$$
Thus, by the construction the disks $$\mathbb
D(\zeta_0,R_{z})\subseteq\mathbb D(\delta_*/|z|)\ ,\ \ \ \mbox{where
$R_{z}\ :=\ \delta_*/|z| - |\zeta_0|\ \geq\ \delta-|\zeta_0|\ >\ 0$\
,\ $z\in\overline{\mathbb D(\tau_*)}\setminus\{ 0\}$ ,}
$$ belong to the domain of definition for the family of the functions $f_{z}(\zeta)$, $z\in\overline{\mathbb
D(\tau_*)}\setminus\{ 0\}$.

 It is clear, see e.g. I.D(8) in
\cite{Ahl}, that $f_{z}(\zeta)$ is a regular homeomorphic solution
of the Beltrami equation with the complex coefficient $\mu_z$ such
that $|\mu_{z}(\zeta)|=|\mu(z\zeta)|$ and
$$
K_{\mu_{z}}(\zeta)\ \leq\ Q_{z}(\zeta)\ :=\ Q(z\zeta)\ \ \ \ \ \
\forall\ \zeta_0\in\mathbb D(\delta)\,,\ \zeta\in\mathbb
D(\zeta_0,R_{z})\ .
$$
Note that the BMO norm of $Q$ as well as its averages over disks are
invariant under linear transformations of variables in $\mathbb C$.
Moreover, the averages $\widetilde Q_{z}(\zeta_0)$ of the function
$Q$ over the disks $\mathbb D(z\zeta_0,\rho_{z})$ forms a continuous
function with respect to the parameter $z\in\overline{\mathbb
D(\tau_*)}\setminus\{ 0\}$ in view of absolute continuity of its
indefinite integrals and it can be extended by continuity to $z=0$
as its (finite !) average over the disk $\mathbb D(\delta_*).$ Since
the closed disk $\overline{\mathbb D(\tau_*)}$ is compact,
$$Q_0\ :=\ \max\limits_{z\in\overline{\mathbb D(\tau_*)}}\,\widetilde Q_{z}(\zeta_0)\ <\ \infty\ .$$

Note also that by Remark 4 $f_{z},\ z\in\overline{\mathbb
D(\tau_*)},$ belongs to the class $\frak{R}^{\Delta}_{Q_{\tau}}$ at
$\zeta_0$ in the punctured disk $\mathbb D(\zeta_0,\delta -
|\zeta_0|)\setminus\{ 0\}$ with $\Delta=1>0$ if $\zeta_0\ne 0$, and
in $D(\zeta_0,\delta - |\zeta_0|)\setminus\{ 1\}$ with
$\Delta=1/\surd 2>1/2$ if $\zeta_0\ne 1$. Hence by Proposition 3 in
any case we obtain the following estimate
\begin{equation}  \label{FMO_eq3.1m}
s(f_{\tau}(\zeta), f_{\tau}(\zeta_0))\ \le\ 64 {\left(\log \frac{2(\delta - |\zeta_0|)}{|\zeta-\zeta_0|}%
\right)^{-\frac{1}{\alpha_0}}}
\end{equation}
for all $z\in\overline{\mathbb D(\tau_*)}$ and $\zeta\in \mathbb
D(\zeta_0,(\delta - |\zeta_0|)/2)$, where $\alpha_0\ =\ 2(Q_0 + 6\|
Q\|_*)$, i.e., the family of the mappings $f_{z}(\zeta),\
z\in\overline{\mathbb D(\tau_*)},$ is equicontinuous at each point
$\zeta_0\in\mathbb D(\delta)$. In view of arbitrariness of
$\delta>1$, the latter is true for all $\zeta_0\in \mathbb C$ at
all.
\end{proof} $\Box$

By the Ascoli theorem, see e.g. 20.4 in \cite{Va}, and Lemma 1 we
obtain the next conclusion.

{\bf Corollary 3.} {\it Let a mapping $f:D\to\mathbb C$ satisfy the
hypotheses of Lemma 1. Then mappings $f_z(\zeta) := f(\zeta \, z) /
f(z)$ form a normal family, i.e., every sequence $f_{z_n}(\zeta)$,
$n=1,2,\ldots$ with $|z_n|\to +0$ as $n\to\infty$ contains a
subsequence $f_{z_{n_k}}(\zeta)$, $k=1,2,\ldots$ that converges with
respect to the spherical metric locally uniformly in $\mathbb C$ as
$k\to\infty$ to a continuous mapping $f_0:\mathbb
C\to\overline{\mathbb C}$ with $f_0(0)=0$ and $f_0(1)=1$. }

\medskip

Further, we are dealing with the so-called approximate solutions of
the Beltrami equations first introduced in the paper \cite{KR}.
Namely, given a domain $D$ in $\mathbb C$, a homeomorphic ACL
(absolutely continuous on lines) solution $f$ of the Beltrami
equation (\ref{eqBeltrami}) in $D$ is called its {\bf approximate
solution} if $f$ is a locally uniform limit in $D$ as $n\to\infty$
of (quasiconformal) homeomorphic ACL solutions $f_n$ of the Beltrami
equations with the complex coefficients
$$\mu_n(z)\ :=\ \begin{cases}\mu(z),& \text{if}\ \mu(z)\leqslant 1-1/n\,,\\ 0,&\text{otherwise\,.}\end{cases}$$
In the paper \cite{KR}, it was established that approximate solution
is unique up to pre-composition with a conformal mapping if
$K_{\mu}\in L^1_{\rm loc}$. Let us give a proof of the following
important fact.

\medskip

{\bf Proposition 4.} {\it Every approximate solution $f$ of Beltrami
equation (\ref{eqBeltrami}) with $K_{\mu}\in L^1_{\rm loc}$ is its
regular homeomorphic solution and, moreover, $f^{\,-1}\in W_{\rm
loc}^{1, 2}$.}

\begin{proof} Indeed, let $f$ be an approximate solution of the Beltrami equation
(\ref{eqBeltrami}) and let $f_n$ be its approximating sequence. Then
first of all $f\in W^{1,1}_{\rm loc}$ by Theorem 3.1 in
\cite{RSY$_8$}.

Let us now prove that $f^{\,-1}\in W_{\rm loc}^{1, 2}.$ Indeed, by
Lemma 3.1 in \cite{RSS} $g_n:=f_{n}^{-1}\rightarrow g:=f^{\,-1}$
uniformly in $\overline{\Bbb C}$ as $n\to\infty$. Note that $f_{n}$
and $g_n\in W_{\rm loc}^{1, 2}$, $n=1,2,\ldots,$ because they are
quasiconformal mappings. Consequently, these homeomorphisms are
locally absolutely continuous, see e.g. Theorem III.6.1 in
\cite{LV}. Observe also that $\mu_n:=(g_n)_{\bar
w}/(g_n)_w=-\mu_{n}\circ g_n,$ see e.g. Section I.C in \cite{Ahl}.
Thus, replacing variables in the integrals, see e.g. Lemma III.2.1
in \cite{LV}), we obtain that
$$\int\limits_{B}|\partial g_n(w)|^{\,2}\,\, dm(w)\ =\ \int\limits_{g_n\left(B\right)}
\frac{dm(z)}{1-|\mu_n(z)|^2}\ \leqslant\ \quad\int\limits_{B^{\,*}}
K_{\mu}(z)\, dm(z)\ <\ \infty$$
for sufficiently large $n,$ where $B$ and $B^{\,*}$ are arbitrary
domains in $\mathbb C$ with compact closures in $f(D)$ and $D$,
respectively, such that $g(\overline{B})\subset B^{\,*}.$ It follows
from the latter that the sequence $g_n$ is bounded in the space
$W^{1,2}(B)$ in each such domain $B$. Hence $f^{\,-1}\in
W^{1,2}_{\rm loc}$, see e.g. Lemma III.3.5 in \cite{Re} or Theorem
4.6.1 in \cite{EG}.

Finally, the latter brings in turn that $g$ has $(N)-$property, see
Theorem III.6.1 in \cite{LV}. Hence $J_f(z)\ne 0$ a.e., see Theorem
1 in \cite{Po}. Thus, $f$ is really a regular solution of the
Beltrami equation (\ref{eqBeltrami}).
\end{proof} $\Box$

\medskip

Note also that Lemma 3.12 in \cite{RSY}, see Lemma 2.12 in the
monograph in \cite{GRSY}, is extended from quasiconformal mappings
to approximate solutions of the Beltrami equation (\ref{eqBeltrami})
im\-me\-dia\-te\-ly by the definition of such solutions.

\medskip

{\bf Proposition 5.} {\it  Let $f:\mathbb D \to \overline{\Bbb C}
\setminus \{ a,b\}, \ a,b \in \overline{\Bbb C}, \ s(a,b) \geq
\delta>0,$ be an  approximate solution of the Beltrami equation
(\ref{eqBeltrami}). Suppose that $s(f(z_1),f(0))\geq \delta$ for
$z_1\in\mathbb D \backslash \{0\}$. Then, for every point $z$ with
$|z| <min(1-|z_1|, |z_1|/2),$
\begin{equation}
s(f(z),f(0)) \geq \psi(|z|)  \label{lasibm_5.24}
\end{equation}
where $\psi$ is a nonnegative strictly increasing function depending
only on $\delta$ and $||K_{\mu}||_1$. }

\medskip

In turn, Propositions 4 and 5 make it possible to prove the
following useful statement.

\medskip

{\bf Proposition 6.} {\it Let $D$ be a domain in $\mathbb{C}$ and
$f_{n}:D\rightarrow\overline{\mathbb C}$ be a sequence of
approximate solutions of the Beltrami equations $\overline{\partial}
f_n=\mu_n\partial f_n$. Suppose that $f_{n}\rightarrow f$ as
$n\to\infty$ locally uniformly in $D$ with respect to the spherical
metric and the norms $\| K_{\mu_n}\|_1$, $n=1,2,\ldots$ are locally
equipotentially bounded. Then either $f$ is constant or it is a
homeomorphism.}

\medskip

\begin{proof} Consider the case when $f$ is not constant in $D$.
Let us first show that then no point in $D$ has a neighborhood of
the constancy for $f$. Indeed, assume that there is at least one
point $z_0\in D$ such that $f(z)\equiv c$ for some
$c\in\overline{\mathbb C}$ in a neighborhood of $z_0$. Note that the
set $\Omega_0$ of such points $z_0$ is open. The set $E_c=\{ z\in D:
s(f(z), c)>0\}$ is also open by continuity of $f$ and not empty if
$f$ is not constant. Thus, there is a point
$z_0\in\partial\Omega_0\cap D$ because $D$ is connected. By
continuity of $f$ we have that $f(z_0)=c$. However, by the
construction there is a point $z_1\in E_c=D\setminus
\overline{\Omega _0}$ such that $|z_0-z_1|<r_0=\mbox{dist}\
(z_0,\partial D)$ and, thus, by the lower estimate of the distance
$s(f(z_0),f(z))$ in Proposition 5 we obtain a contradiction for
$z\in\Omega_0$. Then again by Proposition 5 we obtain that the
mapping $f$ is discrete. Hence $f$ is a homeomorphism by Proposition
3.1 in \cite{RSY$_8$}, see Proposition 2.6 in the monograph
\cite{GRSY}, cf. also Lemma 2.2 and Theorem 2.5 in \cite{GMRV}.
\end{proof} $\Box$

\bigskip

{\bf Corollary 4.} {\it Let a mapping $f:D\to\mathbb C$ satisfy the
hypotheses of Lemma 1 and $f$ be an approximate solution of the
Beltrami equation (\ref{eqBeltrami}) and, moreover,
\begin{equation}\label{eq5.4}
\limsup\limits_{r\to 0}\ \frac{1}{r^2}\int\limits_{|z|<r}\,
|K_{\mu}(z)|\ dm(z)\ <\ \infty\ .
\end{equation}
Then each limit mapping $f_0$ of a sequence $f_{z_n}(\zeta) :=
f(\zeta \, z_n) / f(z_n)$, $z_n\in\mathbb C\setminus\{ 0\}$,
$n=1,2,\ldots$ with $z_n\to 0$ as $n\to\infty$ is a homeomorphism of
$\mathbb C$ into $\mathbb C$.}

\bigskip

\begin{proof} Indeed, $f_{_n}$ are approximate solutions of the Beltrami equations
$\overline\partial f_n=\mu_n\partial f_n$ with
$|\mu_n(\zeta)|=|\mu(z_n\zeta)|$, see e.g. Section I.C in
\cite{Ahl}, and by simple calculations, for all $R>0$,
\begin{equation}\label{eq5.4b}
\overline{\lim\limits_{n\to\infty}}\ \int\limits_{|\zeta|<R}\,
|K_{\mu_n}(\zeta)|\ dm(\zeta)\ =\
R^2\cdot\overline{\lim\limits_{r\to 0}}\
\frac{1}{(R|z_n|)^2}\int\limits_{|z|<R|z_n|}\, |K_{\mu}(z)|\ dm(z)\
<\ \infty
\end{equation}
and, thus,  by Proposition 6 the mapping $f_0$ is a homeomorphism in
$\mathbb C$.


\newpage

Now, let us assume that $f_0(\zeta_0)=\infty$ for some
$\zeta_0\in\mathbb C$. Since $f_n$ are homeomorphisms, there exist
points $\zeta_n\in\mathbb S(\zeta_0,1)$ such that
$s(\zeta_n,\infty)>s(\zeta_0,\infty)$ for all large enough $n$. We
may assume in addition, with no loss of generality, that
$\zeta_n\to\zeta_*\in\mathbb S(\zeta_0,1)$ because the circle
$\mathbb S(\zeta_0,1)$ is a compact set. Then
$f_0(\zeta_*)=\lim\limits_{n\to\infty} f_n(\zeta_n)=\infty$ because
by Lemma 1 the sequence $f_n$ is equicontinuous and, for such
sequences, the pointwise convergence $f_n\to f_0$ is equivalent to
its continuous convergence, see e.g. Theorem 7.1 in \cite{MRSY}.
However, the latter leads to a contradiction because
$\zeta_*\ne\zeta_0$ and by the first part $f_0$ is a homeomorphism.
The obtained contradiction disproves the above assumption and, thus,
really $f_0(\zeta)\ne\infty$ for all $\zeta\in\mathbb C$, i.e.,
$f_0$ is a homeomorphism of $\mathbb C$ into $\mathbb C$.
\end{proof} $\Box$

\section{The main theorems and consequences on asymptotic homogeneity at the origin}

The following theorem shows, in particular, that the Belinskij
conformality still remains to be equi\-va\-lent to the property of
asymptotic homogeneity for regular homeomorphic solutions of the
de\-ge\-ne\-ra\-te Beltrami equations (\ref{eqBeltrami}) if its
dilatation $K_{\mu}$ has a majorant $Q$ in BMO.

\medskip

{\bf Theorem 1.} {\it Let $D$ be a domain in $\mathbb{C}$, $0\in D$,
and let $f:D\rightarrow{\mathbb{C}}$ be a regular homeomorphic
solution of the Beltrami equation with $f(0)=0$ and $K_{\mu}$ have a
majorant $Q\in$ {\rm BMO(D)}. Then the following assertions are
equivalent:

1) $f$ is conformal by Belinskij at the origin,

2) for all $\zeta \in \mathbb{C}$,
\begin{equation}\label{eq_5.5}
    \lim\limits_{\underset{\tau > 0}{\tau \to 0,}} \frac{f(\tau \zeta)}{f(\tau)} =
    \zeta\ ,
\end{equation}

3) for all $\delta>0$, along $z \in \mathbb{C}^{*}:= \mathbb{C}
\setminus \{0\}$ and $z'\in\mathbb C$ with $|z'| \leq \delta |z|$,
\begin{equation}\label{eq_5.6}
 \lim\limits_{z \to 0}
 \left\{ \frac{f(z')}{f(z)} - \frac{z'}{z} \right\}=
 0\ ,
\end{equation}

4) for all $\zeta \in \mathbb{C}$,
\begin{equation}\label{eq_5.7}
   \lim\limits_{\underset{z \in \mathbb{C}^{*}}{z \to 0,}} \frac{f(z \zeta)}{f(z)} =
   \zeta\ ,
\end{equation}

5) the limit in (\ref{eq_5.7}) is uniform in the parameter $\zeta$
on each compact subset of $\mathbb{C}$. }

\medskip

\begin{proof} Let us follow the scheme
 1) $\Rightarrow$ 2) $\Rightarrow$ 3)
$\Rightarrow$ 4) $\Rightarrow$ 5) $\Rightarrow$ 1) and set
$$
f_0(\zeta) = \zeta\ , \ \ \ f_z(\zeta) = f(\zeta \, z) / f(z) \ \ \
\ \ \ \ \mbox{ $\forall\ z\in D\setminus \{0\}$\,, $\ \zeta \in
\mathbb{C}$\,: $z\zeta\in D$\ .}
$$

1) $\Rightarrow$ 2). Immediately the definition of the conformality
by Belinskij yields the convergence $f_{\tau} (\zeta) \to
f_0(\zeta)$ as $\tau \to 0$ along $\tau >0 $ for every fixed $\zeta
\in \mathbb{C}$, i.e., just (\ref{eq_5.5}).

2) $\Rightarrow$ 3). In view of Lemma 1, the pointwise convergence
in (\ref{eq_5.5}) for each $\zeta\in\mathbb C$ implies the uniform
convergence there on compact sets in $\mathbb C$, see e.g. Theorem
7.1 in \cite{MRSY}. To obtain on this basis the implication 2)
$\Rightarrow$ 3), let us note the identities
$$
f_z(\zeta)\ =\ \frac{f_{|z|} (\zeta z / |z|)}{f_{|z|} (z / |z|)}\ =\
\frac{f(z')}{f(z)} \ \ \ \ \ \ \ \forall\ \zeta = \frac{z'}{z} \in
\mathbb{C}\, ,\ z\in\mathbb{C}^{*} := \mathbb{C} \setminus \{0\}\ .
$$
Hence to prove (\ref{eq_5.6}) it is sufficient to show that $f_z
(\zeta) - f_0 (\zeta) \to 0$ as $z \to 0, \ z\in \mathbb{C}^*$
uniformly with respect to the parameter $\zeta$ in the closed disks
$\mathbb D_{\delta} := \{ \zeta \in \mathbb{C} : |\zeta| \leq \delta
\}$, $\delta > 0$.

Indeed, let us assume the inverse. Then there is a number
$\varepsilon > 0$ and consequences $\zeta_n \in D_{\delta}, \ z_n
\to 0, \ z_n \in \mathbb{C}^{*},$ such that $|g_n (\zeta_n) -
\zeta_n | \geq \varepsilon$, where $g_n (\zeta) = f_{z_n} (\zeta), \
\zeta \in \mathbb{C}$. Since the closed disk $\mathbb D_{\delta}$
and the unit circle $\partial \mathbb D_1$ are compact sets, then
with no loss of generality we may in addition to assume that
$\zeta_n \to \zeta_0 \in\mathbb D_{\delta}$ and $\eta_n = z_n /
|z_n| \to \eta_0 \in \partial\mathbb D_1$ as $n \to \infty$.

Let us denote by $\varphi_n (\zeta)$ the mappings $f_{|z_n|}
(\zeta), \ \zeta \in \mathbb{C}, \ n=1, \, 2, \, \dots.$ Then
$\varphi_n(\zeta) \to \zeta$ as $n \to \infty$ uniformly on $\mathbb
D_{\delta} \cup \, \partial\mathbb D_1$ and $g_n (\zeta) = \varphi_n
(\eta_n \zeta) / \varphi_n(\eta_n)$. Consequently, $g_n (\zeta) \to
\zeta$ as $n \to \infty$ uniformly on $\mathbb D_{\delta}$. Hence
$g_n (\zeta_n) \to \zeta_0$ as $n \to \infty$ because the uniform
convergence of continuous mappings on compact sets implies the
so-called continuous convergence, see e.g. Remark 7.1 in
\cite{MRSY}. Thus, the obtained contradiction disproves the above
assumption.

3) $\Rightarrow$ 4). Setting in (\ref{eq_5.6}) $z' = z \zeta$ and
$\delta = |\zeta|$, we immediately obtain (\ref{eq_5.7}).

4) $\Rightarrow$ 5). The limit relation (\ref{eq_5.7}) means in the
other words that $f_z(\zeta) \to f_0(\zeta)$ as $z \to 0$ along $\ z
\in \mathbb{C}^*$ pointwise in $\mathbb{C}$.  In view of Lemma 1,
the latter implies the locally uniform convergence $f_z \to f_0$ as
$z \to 0$ in $\mathbb{C}$, see again Theorem 7.1 in \cite{MRSY}.

5) $\Rightarrow$ 1). From (\ref{eq_5.7}) for $z = \rho > 0, \ \zeta
= e^{i \vartheta}, \ \vartheta \in \mathbb{R},$ and $w = \zeta z =
\rho \, e^{i \vartheta}$ we obtain that $f(w) = f(\rho) (\zeta +
\alpha (\rho)),$ where $\alpha(\rho) \to 0$ as $\rho \to 0$.
Consequently,
$$
f(w) = A(\rho) (w + o(\rho))\ ,
$$
where $ A(\rho) = {f(\rho)}/{\rho} $ and $o(\rho) / \rho \to 0$ as
$\rho \to 0$. Moreover, by (\ref{eq_5.7}) with $z= \rho
> 0$ and $\zeta = t > 0$ we have that $A$ satisfies the condition
$$
\lim\limits_{\rho \to 0} \frac{A(t \rho)}{A(\rho)} = 1 \quad
    \forall\,t>0\ ,
$$
i.e., $f$ is conformal by Belinskij at the origin. \end{proof}
$\Box$

\medskip

The following result is fundamental for further study of asymptotic
homogeneity because it facilitates considerably the verification of
(\ref{eq_5.7}) and at the same time reveals the nature of the
notion. Let $Z$ be an arbitrary set in the complex plane $\mathbb
C$, $0\notin Z,$ with the origin as its accumulation point. Further,
we use the following characteristic of its sparseness:
\begin{equation}\label{eq_5.8}
\mathbb S_Z(\rho)\ :=\ \frac{\inf_{{z \in Z},{|z|\geq
\rho,}}|z|}{\sup_{{z \in Z},{|z|\leq \rho,}}|z|} \ \ \ \ \  \forall\
\rho\ > \ 0\ .
\end{equation}

{\bf Theorem 2.} {\it Let $f$ satisfy the hypotheses of Theorem 1.
Suppose that
\begin{equation}\label{eq_5.9}
\limsup\limits_{\rho\to 0}\ \mathbb S_Z(\rho)\ <\ \infty\
\end{equation}
and
\begin{equation}\label{eq_5.10}
   \lim\limits_{\underset{z \in Z}{z \to 0,}} \frac{f(z \zeta)}{f(z)} =
   \zeta\ \ \ \ \ \ \forall\ \zeta\in\mathbb C\ .
\end{equation}
Then $f$ is asymptotically homogeneous at the origin.}

\medskip

{\bf Remark 4.} For Theorem 2 to be true, the condition
(\ref{eq_5.9}) on the extent of possible sparseness of $Z$ is not
only sufficient but also necessary as Proposition 2.1 in \cite{GR}
in the case $Q\in L^{\infty}\subset$ BMO shows. In particular, any
continuous path to the origin or a discrete set, say $1/n , n = 1 ,
2 , \ldots$, can be taken as the set $Z$ in Theorem 2. For instance,
the conclusion of Theorem 2 is also true if $Z$ has at least one
point on each circle $|z|=\rho$ for all small enough $\rho>0$.

\medskip

\begin{proof}
Indeed, by (\ref{eq_5.10}) we have that, for functions $f_z(\zeta)
:=f(z \zeta)/f(z)$, pointwise
\begin{equation}\label{eq_5.11}
   \lim\limits_{\underset{z \in Z}{z \to 0,}} f_z(\zeta)\ =\
   \zeta\ \ \ \ \ \ \forall\ \zeta\in\mathbb C
\end{equation}
and, by Theorem 7.1 in \cite{MRSY} and Lemma 1, the limit in
(\ref{eq_5.11}) is locally uniform in $\zeta\in\mathbb C$.

Let us assume that (\ref{eq_5.7}) does not hold for $f$, in the
other words, there exist $\zeta\in\mathbb C$, $\varepsilon > 0$ and
a sequence $z_n\in\mathbb C^*$, $n=1,2,\ldots$ such that $z_n\to 0$
as $n\to\infty$ and
\begin{equation}\label{eq_5.12}
 |\ f_{z_n}(\zeta)\ -\ \zeta\ |\ >\ \varepsilon\ .
\end{equation}
On the other hand, by (\ref{eq_5.9}) there is a sequence $z^*_n\in
Z$ such that
$$
0\ <\ \delta\ \leq\ |\tau_n|\ \leq\ 1\ <\ \infty
$$
for all large enough $n=1,2,\ldots $, where
$$
\tau_n\ =\ \frac{z_n}{z_n^*}\ ,\ \ \ \delta\ =\
1/2\limsup\limits_{\rho\to 0}\, \mathbb S_Z(\rho)\ .
$$
With no loss of generality, we may assume in addition that
$\tau_n\to\tau_0$ with $\delta\leq|\tau_0|\leq 1$ as $n\to\infty$
because the closed ring $R:=\{z\in\mathbb C:\delta\leq|z|\leq 1\}$
is a compact set. Note also that
$$
f_{z_n}(\zeta)\ =\ \frac{f_{z^*_n}(\zeta\tau_n)}{f_{z^*_n}(\tau_n)}\
.
$$
Thus, $f_{z^*_n}(\zeta\tau_n)\sim \zeta\tau_0$ and
$f_{z^*_n}(\tau_n)\sim \tau_0$ as $n\to\infty$ because the uniform
convergence in (\ref{eq_5.11}) with respect to $\zeta$ over any
compact set implies the so-called continuous convergence, see e.g.
Remark 7.1 in \cite{MRSY}. Consequently, $f_{z_n}(\zeta)\sim\zeta$
as $n\to\infty$ because $\tau_0\ne0$. However, the latter
contradicts (\ref{eq_5.12}). The obtained contradiction disproves
the above assumption and the conclusion of the theorem is true.
\end{proof} $\Box$

\medskip

Now, recall that the abstract spaces $\frak F$ in which convergence
is a primary notion were first considered by Frechet in his thesis
in 1906. Later on, Uryson introduced the third axiom in these
spaces: if a compact sequence $f_n\in\frak F$ has its unique
accumulation point $f\in\frak F$, then $\lim\limits_{n\to\infty} f_n
= f$, see e.g. \cite{Ku}, Chapter 2, § 20,1-II. Recall that
$f_n\in\frak F$, $n=1,2,\ldots$ is called a {\bf compact sequence}
if each its subsequence contains a converging subsequence and,
moreover, $f\in\frak F$ is said to be an {\bf accumulation point} of
the sequence $f_n\in\frak F$ if $f$ is a limit of some its
subsequence. It is customary to call such spaces {\bf ${\mathfrak
L}^*-$spaces}.

{\bf Remark 5.} {\it In particular, any convergence generated by a
metric satisfies Uryson's axiom}, see e.g. \cite{Ku}, Chapter 2, §
21, II. However, the well-known convergence almost everywhere of
measurable functions yields a counter-example to Uryson's axiom: any
sequence converging in measure is compact with respect to
convergence almost everywhere, but not every such sequence converges
almost everywhere. Later on, we apply the convergence generated by
the uniform convergence of continuous functions, generated as known
by the uniform norm.

\medskip

To prove the corresponding sufficient criteria for the asymptotic
homogeneity at the origin for solutions of degenerate Beltrami
equations, we need also the following general lemma.

\medskip

{\bf Lemma 2.} {\it Let $D$ be a bounded domain in $\mathbb{C}$ and
$f_{n}:D\rightarrow{\mathbb C}$, $n=1,2,\ldots$ be a sequence of
$W^{1,1}$ solutions of the Beltrami equations $\overline{\partial}
f_n=\mu_n\partial f_n$. Suppose that $f_{n}\rightarrow f$ as
$n\to\infty$ in $L^1$ and the norms $\| \overline{\partial} f_n\|_1$
and $\| \overline{\partial} f_n\|_1$ are equipotentially bounded.
Then $f\in W^{1,1}$ and $\partial f_{n}$ and $\overline{\partial}
f_{n}$ converge weakly in $L^1$ to $\partial f$ and
$\overline{\partial} f$, respectively. Moreover, if
$\mu_{n}\rightarrow\mu$ a.e. or in measure as $n\to\infty$, then
$\overline{\partial} f=\mu\partial f$ a.e.}

\medskip

\begin{proof} The first part of conclusions follow from Lemma
III.3.5 in \cite{Re}. Let us prove the latter of these conclusions.
Namely, assuming that $\mu_{n}(z)\rightarrow\mu(z)$ a.e. as
$n\to\infty$ and, setting $$\zeta(z)\ =\ \overline{\partial}f(z) -
\mu(z)\cdot\partial f(z)\ ,$$ let us show that $\zeta(z)=0$. Indeed,
since $\overline{\partial} f_n(z)-\mu_n(z)\partial f_n(z)=0$, by the
triangle inequality
$$\int\limits_{D}|\zeta(z)|\,dm(z)\ \leq\ I_1(n)\ +\ I_2(n)\ +\ I_3(n)\ ,$$
where
$$I_1(n)\ :=\ \int\limits_{D} |\overline{\partial}f(z) - \overline{\partial}f_n(z)|\,dm(z)\ ,$$
$$I_2(n)\ :=\ \int\limits_{D}|\mu(z)|\cdot|\partial f(z) - \partial
f_n(z)|\,dm(z)\ ,$$
$$I_3(n)\ :=\ \int\limits_{D}|\mu(z)-\mu_n(z)|\cdot|\partial f_n(z)|\,dm(z)\ .$$
By the first part of conclusions, with no loss of generality, assume
that $|\overline{\partial}f(z)-\overline{\partial}f_n(z)|\to 0$ and
$|\partial f(z)-\partial f_n(z)|\to 0$ as $n\to\infty$ weakly in
$L^1$, see Corollary IV.8.10 in \cite{DS}. Thus, $I_1(n)\to 0$ and
$I_2(n)\to 0$ as $n\to\infty$ because the dual space of $L^1$ is
naturally isometric to $L^{\infty}$, see e.g. Theorem IV.8.5 in
\cite{DS}.

Moreover, by Corollary IV.8.11 in \cite{DS}, for each $\varepsilon>
0$, there is $\delta>0$ such that over every measurable set $E$ in
$D$ with $|E|<\delta$
\begin{equation}\label{dop2e3.1.10}
\int\limits_{E}|\
\partial f_n(z)|\,dm(z)<\varepsilon,\quad n=1,2,\ldots\ .
\end{equation} Further, by the Egoroff theorem, see e.g. III.6.12 in \cite{DS},
$\mu_n(z)\to\mu(z)$ as $n\to\infty$ uniformly on some set $S$ in $D$
with $|E| < \delta$ where $E = D\setminus S$. Hence  $|\mu_n(z) -
\mu(z)| < \varepsilon$ on $S$ and
$$
I_3(n) \leq \varepsilon\int\limits_{S}\ |\ \partial f_n(z) |\
\,dm(z)+2\int\limits_{E}|\ \partial f_n(z)|\ \,dm(z)\ \leq\
\varepsilon(\|\partial f_n(z)\|_1\ +\ 2)
$$
for large enough $n$, i.e., $I_3(n)\to 0$ because $\varepsilon >0$
is arbitrary. Thus, really $\zeta = 0$ a.e.
\end{proof} $\Box$

\medskip

{\bf Theorem 3.} {\it Let $D$ be a domain in $\mathbb{C}$, $0\in D$,
$f:D\rightarrow{\mathbb{C}}$, $f(0)=0$, be an approximate solution
of the Beltrami equation (\ref{eqBeltrami}) and $K_{\mu}$ have a
majorant $Q\in$ {\rm BMO(D)}. Suppose that
\begin{equation}\label{eq5.4}
\limsup\limits_{r\to 0}\ \frac{1}{\pi r^2}\int\limits_{|z|<r}\,
K_{\mu}(z)\ dm(z)\ <\ \infty\ ,
\end{equation}
and
\begin{equation}\label{eq5.5}
\lim\limits_{r\to 0}\ \frac{1}{\pi r^2}\int\limits_{|z|<r}\,
|\mu(z)|\ dm(z)\ =\ 0\ .
\end{equation}
Then $f$ is asymptotically homogeneous at the origin. }

\medskip

\begin{proof} By Theorem 2 with $Z:=\{ 2^{-n}\}^{\infty}_{n=N}$, where $2^{-N}<{\rm dist}(0,\partial D)$, it is sufficient to show
that
$$
\lim\limits_{n\to \infty}\ f_{2^{-n}}(\zeta)\ =\ \zeta\ \ \ \ \
\forall\ \zeta\in\mathbb C\ ,\ \ \ f_{2^{-n}}(\zeta)\ :=\
\frac{f(2^{-n}\zeta)}{f(2^{-n})}\ .
$$
By Corollary 3 the sequence $f_{2^{-n}}(\zeta)$ is compact with
respect to locally uniform convergence in $\mathbb C$ and by Remark
5 it remains to prove that each its converging subsequence
$f^*_k=f_{n_k}$ with $n_k\to\infty$ as $k\to\infty$ has the identity
mapping of the complex plane $\mathbb C$ as its limit $f_0$.

Indeed, the mappings $f^*_k$ are approximate solutions of Beltrami
equations $\overline{\partial}f^*_k=\mu^*_k\cdot\partial f^*_k$ with
$|\mu^*_k(\zeta)|=|\mu(2^{-n_k}\zeta)|$, see e.g. calculations of
Section I.C in \cite{Ahl}. Since such solutions are regular by
Proposition 4, we have by the calculations that
$$|\ \overline{\partial}f^*_k|\leq|\ \partial f^*_k|\leq|\ \partial
f^*_k|+|\ \overline{\partial}f^*_k| \leq K_{\mu^*_k}^{1/2}
J_{f^*_k}^{1/2}\ \ \ \  \ \ \  \textrm{a.e.}\ ,\ \ \ k=1,2,\ldots$$
where
$$K_{\mu^*_k}(\zeta)=K_{\mu}(2^{-n_k}\zeta)\ ,\ J_{f^*_k}(\zeta)=|\partial
f^*_k(\zeta)|^2-|\overline{\partial}f^*_k(\zeta)|^2=J_{f_{n_k}}(\zeta)=J_{f}(2^{-n_k}\zeta)/|f(2^{-n_k})|^2
.$$ Consequently, by the H\"older inequality for integrals, see e.g.
Theorem 189 in \cite{HLP}, and Lemma III.3.3 in \cite{LV}, we obtain
that
$$\| \partial f^*_{k}\|_1(\mathbb D_l)\ \leq\ \| K_{\mu^*_k}\|^{\frac{1}{2}}_1(\mathbb D_l) \cdot\left\vert
f^*_{k}(\mathbb D_l) \right\vert ^{\frac{1}{2}}\ \ \ \ \ \ \forall\
l=1,2,\ldots\ ,\ \mathbb D_l:=\mathbb D(2^l)\ .$$

Now, by the condition (\ref{eq5.4}) and simple calculations, for
each fixed $l=1,2,\ldots\ ,$
$$
\overline{\lim\limits_{k\to\infty}}\ \| K_{\mu^*_k}\|_1(\mathbb
D_l)\ =\ 2^{2l}\cdot\overline{\lim\limits_{k\to \infty}}\
\frac{1}{(2^l2^{-n_k})^2}\int\limits_{|z|<2^l2^{-n_k}}\,
|K_{\mu}(z)|\ dm(z)\ <\ \infty\ .
$$
Next, choosing $\zeta_k$ in $\mathbb S_l:=\{\zeta\in\mathbb
C:|\zeta|=2^l\}$ with
$|f(2^{-n_k}\zeta_k)|=\max\limits_{\zeta\in\mathbb
S_l}|f(2^{-n_k}\zeta)|$, we see that
$$
\left\vert f^*_{k}(\mathbb D_l) \right\vert = \left\vert
f_{n_k}(\mathbb D_l) \right\vert = \frac{\left\vert f(\mathbb
D(2^{l-n_k}))
\right\vert}{|f(2^{-n_k}\zeta_k)|^2}\cdot\frac{|f(2^{-n_k}\zeta_k)|^2}{|f(2^{-n_k})|^2}\
\leq\ \pi |f_{2^{-n_k}}(\zeta_k)|^2\ .
$$
With no loss of generality, we may assume that $\zeta_k\to\zeta_0\in
\mathbb S_l$ as $k\to\infty$ because the circle $\mathbb S_l$ is a
compact set. Then $f_{2^{-n_k}}(\zeta_k)\to f_0(\zeta_0)$ because
the uniform convergence implies the so-called continuous
convergence, see e.g. Remark 7.1 in \cite{MRSY}. However,
$f_0(\zeta_0)\ne\infty$, see Corollary 4.

Thus, the norms  of ${\partial}f^*_k$ and $\overline{\partial}f^*_k$
are locally equipotentially bounded in $L^1$. Then $f_0$ is
$W^{1,1}_{\rm loc}$ solution of the Beltrami equation with
$\mu\equiv 0$ in $\mathbb C$ by Lemma 2 in view of (\ref{eq5.5}).
Moreover, $f_0$ is a homeomorphism of $\mathbb C$ into $\mathbb C$
by Corollary 4. Hence $f_0$ is a conformal mapping of $\mathbb C$
into $\mathbb C$, see e.g. Corollary II.B.1 in \cite{Ahl}. Hence
$f_0(\zeta)$ is a linear function $a+b\zeta$, see e.g. Theorem
2.31.1 in \cite{LS}. In addition, by the construction $f_0(0)=0$ and
$f_0(1)=1$. Thus, $f_0(\zeta)\equiv\zeta$ in the whole complex plane
$\mathbb C$ and the proof is thereby complete.
\end{proof} $\Box$

\medskip

{\bf Remark 6.} Note that, in particular, both conditions
(\ref{eq5.4}) and (\ref{eq5.5}) follow from the only one stronger
condition
\begin{equation}\label{eq5.6}
\lim\limits_{r\to 0}\ \frac{1}{\pi r^2}\int\limits_{|z|<r}\,
K_{\mu}(z)\ dm(z)\ =\ 1
\end{equation}
because
\begin{equation}\label{eq5.7}
|\mu(z)|\ \leq\ \frac{|\mu(z)|}{1-|\mu(z)|}\ =\
\frac{K_{\mu}(z)-1}{2}\ .
\end{equation}

\medskip

Combinig Theorems 1 and 3, see also Proposition 4, we obtain the
following conclusions.

\medskip

{\bf Corollary 5.} {\it Under hypotheses of Theorem 3, $f$ is
conformal by Lavrent'iev at the origin, i.e., $f$ preserves
infinitesimal circles centered at the origin:
\begin{equation}\label{eq_4.8d}
    {\ \lim\limits_{r \to 0}\ \frac{\max\limits_{|z| =r} |f(z)|} {\min\limits_{|z| =r} |f(z)|}\ =\
    1\ }\ ,
\end{equation}
asymptotically preserves angles, i.e.,
\begin{equation}\label{eq_4.9d}
    {\ \lim\limits_{z \to 0}\ \arg\, \left[\frac{f(z \zeta)}{f(z)} \right] = \arg
    \zeta\ }\ \ \ \ \ \forall\ \zeta \in \mathbb{C}\ ,\ |\zeta|\ =\
    1\ ,
\end{equation}
and asymptotically preserves the moduli of infinitesimal rings,
i.e.,
\begin{equation}\label{eq_4.10d}
    {\ \lim\limits_{z \to 0}\ \frac{|f(z \,\zeta)|}{|f(z)|}\ =\
    |\zeta|\ }\ \ \ \ \ \forall\ \zeta \in \mathbb{C}^*\ :=\ \mathbb C\setminus\{ 0\}\ .
\end{equation}}


{\bf Corollary 6.} {\it Under hypotheses of Theorem 3, for all
$\delta>0$, along $z \in \mathbb{C}^{*}:= \mathbb{C} \setminus
\{0\}$ and $z'\in\mathbb C$ with $|z'| \leq \delta |z|$,
\begin{equation}\label{eq_5.6m}
 \lim\limits_{z \to 0}
 \left\{ \frac{|f(z')|}{|f(z)|} - \frac{|z'|}{|z|} \right\}=
 0\ .
\end{equation}}

Moreover, by the theorem of Stolz (1885) and Cesaro (1888), see e.g.
Problem 70 in \cite{PS}, we derive from Corollary 6 the next
assertion on logarithms.


{\bf Corollary 7.} {\it Under hypotheses of Theorem 3,
\begin{equation}\label{eq_4.43}
\lim\limits_{\underset{z \in \mathbb{C}^{*}}{z \to 0}} \frac{\ln
|f(z)|}{\ln |z|}\ =\ 1\ .
\end{equation}}

\begin{proof} For brevity, let us introduce designations $t_n = - \ln |z_n|$, $\tau_n = -\ln
|f(z_n)|$ and assume that (\ref{eq_4.43}) does not hold, i.e., there
exist $\varepsilon> 0$ and a sequence $z_n \to 0$ such that
\begin{equation}\label{eq_4.45}
    \left|\frac{\tau_n}{t_n} - 1 \right| \geq \varepsilon\ \ \ \ \
    \forall\ n=1, \, 2, \, \dots\ .
\end{equation}
Passing, if necessary, to a subsequence, we can consider that $t_n -
t_{n-1} \geq 1$ for all $n=1, \, 2, \, \dots$. Then, we can achieve
that $t_n - t_{n-1} < 2$, by inserting, if necessary, the mean
arithmetic values between neighboring terms of the subsequence $t_n,
\ n=1, \, 2, \, \dots$. In this case, inequality (\ref{eq_4.45})
holds for the infinite number of terms of the subsequence.

Thus, the sequence $\rho_n = |z_n| = e^{-t_n}$ satisfies the
inequalities $e^{-2} < \rho_n   / \rho_{n-1} \leq e^{-1}$. Relations
(\ref{eq_5.6m}) implies that $\exp (\tau_{n-1} - \tau_n) = \exp
(t_{n-1} - t_n) + \alpha_n,$ where $\alpha_n \to 0$ as $n \to
\infty$, or, in the other form, $\exp (\tau_{n-1} - \tau_n) = (1 +
\beta_n) \exp (t_{n-1} - t_n)$ with $\beta_n \to 0$ as $n \to
\infty$. The latter gives that $(\tau_{n} - \tau_{n-1}) = (t_{n} -
t_{n-1}) + \gamma_n$ with $\gamma_n \to 0$ as $n \to \infty$ and,
since $t_n - t_{n-1} \geq 1,$ we have that $\ (\tau_n - \tau_{n-1})
/ (t_n - t_{n-1}) = 1+ \delta_n,$ where $\delta_n \to 0$ as $n \to
\infty$. By the Stolz theorem, then we conclude that $\tau_n / t_n
\to 1$ in contradiction with (\ref{eq_4.45}). This contradiction
disproves the above assumption, i.e., (\ref{eq_4.43}) is true.
\end{proof} $\Box$


{\bf Theorem 4.} {\it Let $D$ be a domain in $\mathbb{C}$ and let
$f:D\rightarrow{\mathbb{C}}$ be an approximate solution of the
Beltrami equation (\ref{eqBeltrami}), $K_{\mu}$ have a majorant
$Q\in$ {\rm BMO(D)} and  at a point $z_0\in D$
\begin{equation}\label{eq5.4Bel}
\limsup\limits_{r\to 0}\ \frac{1}{\pi r^2}\int\limits_{|z-z_0|<r}\,
K_{\mu}(z)\ dm(z)\ <\ \infty\ .
\end{equation}
Suppose that $\mu(z)$ is approximately continuous at $z_0$. Then the
mapping $f$ is differentiable by Belinskij at this point with $\mu_0
= \mu(z_0)$.}

\medskip

\begin{proof}
First of all, $|\mu(z_0)|<1$ because by the hypotheses $K_{\mu}\in
L^1_{\rm loc}$ and $\mu(z)$ is approximately continuous at $z_0$.
Note also that $f$ is differentiable by Belinskij with $\mu_0 =
\mu(z_0)$ at $z_0$ if and only if $g := h \circ \varphi^{-1}$ is
conformal by Belinskij at zero, where $ h(z) = f(z_0 + z) - f(z_0)$
and $\varphi(z) = z + \mu_0 \overline{z}.$ It is evident that
$\mu_h(z)=\mu(z+z_0)$ and $K_{\mu_h}=K_{\mu}(z+z_0)$ and by
elementary calculations, see e.g. Section I.C(6) in \cite{Ahl},
$\mu_g\circ\varphi=(\mu_h-\mu_0)/(1-\overline{\mu_0}\mu_h)$ and
$K_{\mu_g}\leq K_0\cdot K_{\mu_h}\circ\varphi^{-1}\leq K_0Q_0$,
where $K_0=(1+|\mu_0|)/(1-|\mu_0|)$ and
$Q_0(w)=Q(z_0+\varphi^{-1}(w))$ belongs to BMO in $D_0:=\varphi(D)$
because $\varphi$ and $\varphi^{-1}$ are $K_0-$quasiconformal
mappings, see the paper \cite{Rei} and the monograph \cite{RR}.
Thus, Theorem 4 follows from Theorem 3.
\end{proof} $\Box$

\section{On homeomorphic solutions in extended complex plane}

Here we start from establishing a series of criteria for existence
of approximate solutions $f:\mathbb C\to\mathbb C$ to the degenerate
Beltrami equations in the whole complex plane $\mathbb C$ with the
normalization $f(0)=0$, $f(1)=1$ and $f(\infty)=\infty$.

It is easy to give examples of locally quasiconformal mappings of
$\mathbb C$ onto the unit disk $\mathbb D$, consequently, there
exist locally uniform elliptic Beltrami equations with no such
solutions. Hence, compared with our previous articles, the main goal
here is to find the corresponding additional conditions on
dilatation quotients of the Beltrami equations at infinity.

\bigskip

{\bf Lemma 3.} {\it Let a function $\mu : \mathbb C\to{\Bbb C}$ be
measurable with $|\mu (z)| < 1$  a.e., $K_{\mu}\in L^1_{\rm
loc}(\mathbb C).$ Suppose that, for every $z_0\in\overline{\mathbb
C}$, there exist $\varepsilon_0=\varepsilon(z_0)>0$ and a family of
measurable functions
${\psi}_{z_0,{\varepsilon}}:(0,\infty)\to(0,\infty)$ such that
\begin{equation}\label{eqI}
I_{z_0}({\varepsilon})\ \colon =\
\int\limits_{{\varepsilon}}^{{\varepsilon}_0}{\psi}_{z_0,{\varepsilon}}(t)\
dt\ <\ \infty\ \ \ \ \ \ \forall\
{\varepsilon}\in(0,{\varepsilon}_0)\end{equation} and
\begin{equation}\label{eqT}
\int\limits_{{\varepsilon}<|z-z_0|<{\varepsilon}_0}\
K^T_{{\mu}}(z,z_0)\cdot{\psi}^2_{z_0,{\varepsilon}}(|z-z_0|)\ dm(z)\
=\ o(I^2_{z_0}({\varepsilon})) \ \ \  \ \ \ \hbox{ as
${\varepsilon}\to 0$}\ \ \ \forall\ z_0\in \mathbb C\
\end{equation}
and, moreover,
\begin{equation}\label{eqINFTY}
\int\limits_{{\varepsilon}<|\zeta|<{\varepsilon}_{\infty}}\
K^T_{{\mu}}(\zeta,\infty)\cdot{\psi}^2_{\infty,{\varepsilon}}(|\zeta|)\
\frac{dm(\zeta)}{|\zeta|^4}\ =\ o(I^2_{\infty}({\varepsilon})) \ \ \
\ \ \ \hbox{ as ${\varepsilon}\to 0$}\ ,
\end{equation}
where $K^T_{{\mu}}(\zeta,\infty):=K^T_{{\mu}}({1}/{{\zeta}},0)$.

Then the Beltrami equation (\ref{eqBeltrami}) has an approximate
homeomorphic solution $f$ in $\mathbb C$ with the
nor\-ma\-li\-za\-tion $f(0)=0$, $f(1)=1$ and $f(\infty)=\infty$.}

\bigskip

{\bf Remark 7.} After the replacements of variables $\zeta
\longmapsto z:=1/\zeta$, $\varepsilon\longmapsto R:=1/\varepsilon$,
$\varepsilon_{\infty}\longmapsto R_0:=1/\varepsilon_{\infty}$ and
functions
$\psi_{\infty,\varepsilon}(t)\longmapsto\psi_R(t):=\psi_{\infty,1/R}(1/t)$,
the condition (\ref{eqINFTY}) can be rewritten in the more
convenient form:
\begin{equation}\label{eqINFTYre}
\int\limits_{R_0<|z|<R}\ K^T_{{\mu}}(z,0)\ {\psi}^2_{R}(|z|)\
\frac{dm(z)}{|z|^4}\ =\ o(I^2(R)) \ \ \ \ \ \ \hbox{ as $R\to\infty
$}\ ,
\end{equation}
with the family of measurable functions
${\psi}_{R}:(0,\infty)\to(0,\infty)$ such that
\begin{equation}\label{eqIre}
I(R)\ \colon =\ \int\limits_{R_0}^{R}{\psi}_{R}(t)\ \frac{dt}{t^2}\
<\ \infty\ \ \ \ \ \ \forall\ R\in(R_0,\infty)\ .\end{equation}

Before to come to the proof of Lemma 3, let us recall that a {\bf
condenser} in ${\mathbb C}$ is a domain ${\cal R}$ in ${\mathbb C}$
whose complement in $\overline{\mathbb C}$ is the union of two
distinguished disjoint compact sets $C_1$ and $C_2$. For
convenience, it is written ${\cal R}={\cal R}(C_1,C_2)$. A {\bf
ring} in ${\mathbb C}$ is a condenser  ${\cal R}={\cal R}(C_1,C_2)$
with connected  $C_1$ and $C_2$ that are called the {\bf
complementary components} of ${\cal R}$. It is known that the
(conformal) capacity of a ring   ${\cal R}={\cal R}(C_1,C_2)$ in
$\mathbb C$ is equal to the (conformal) modulus of all paths in
${\cal R}$ connecting $C_1$ and $C_2$, see e.g. Theorem A.8 in
\cite{MRSY}.

\bigskip

\begin{proof}
By the first item of the proof of Lemma 3 in \cite{RSY$_4$} the
Beltrami equation (\ref{eqBeltrami}) has  under the conditions
(\ref{eqT}) an approximate homeomorphic solution $f$ in $\mathbb C$
with $f(0)=0$ and $f(1)=1$. Moreover, by Lemma 3 in \cite{RSY$_4$}
we may also assume that $f$ is a ring $Q-$homeomorphism with
$Q(z)=K^T_{\mu}(z,0)$ at the origin, i.e., for every ring
$A=A(r_1,r_2):=\{ z\in\mathbb C: r_1<|z|<r_2\}$, we have the
estimate of the capacity $C_f(r_1,r_2)$ of its image under the
mapping $f$:
$$
C_f(r_1,r_2)\ \leq\ \int\limits_{A(r_1,r_2)}K^T_{\mu}(z,0)\ dm(z)\ \
\ \ \ \ \forall\ r_1,r_2:\ 0<r_1<r_2<\infty\ .
$$

Let us consider the mapping $F(z):=1/f(1/z)$ in $\mathbb
C_*:=\overline{\mathbb C}\setminus\{ 0\}$. Note that
$F(\infty)=\infty$ because $f(0)=0$. Since the capacity is invariant
under conformal mappings, we have by the change of variables
$z\longmapsto\zeta:=1/z$ as well as
$r_1\longmapsto\varepsilon_2:=1/r_1$ and
$r_2\longmapsto\varepsilon_1:=1/r_2$ that
$$
C_F(\varepsilon_1,\varepsilon_2)\ \leq\
\int\limits_{A(\varepsilon_1,\varepsilon_2)}K^T_{\mu}(1/{\zeta},0)\
\frac{dm(\zeta)}{|\zeta|^4}\ \ \ \ \ \ \forall\
\varepsilon_1,\varepsilon_2:\ 0<\varepsilon_1<\varepsilon_2<\infty\
,
$$
i.e., $F$ is a ring $\tilde Q-$homeomorphism at the origin with
$\tilde Q(\zeta):=K^T_{\mu}(1/{\zeta},0)/|\zeta|^4$. Thus, in view
of the condition (\ref{eqINFTY}), we obtain by Lemma 6.5 in
\cite{GRSY} that $F$ has a continuous extension to the origin. Let
us assume that $c:=\lim\limits_{\zeta\to 0}F(\zeta)\ne 0$.

However, $\overline{\mathbb C}$ is homeomorphic to the sphere
$\mathbb S^2$ by stereographic projection and hence by the Brouwer
theorem in $\mathbb S^2$ on the invariance of domain the set $C_*:=
F({\mathbb C_*})$ is open in $\overline{\mathbb C}$, see e.g.
Theorem 4.8.16 in \cite{Sp}. Consequently, $c\notin C_*$ because $F$
is a homeomorphism. Then the extended mapping $\tilde F$ is a
homeomorphism of $\overline{\mathbb C}$ into $\mathbb C_*$ because
$f\ne \infty$ in $\mathbb C$. Thus, again by the Brouwer theorem,
the set $C:=\tilde F(\overline{\mathbb C})$ is open in
$\overline{\mathbb C}$ and $0\in\overline{\mathbb C}\setminus
C\ne\emptyset$. On the other hand, the set $C$ is compact as a
continuous image of the compact space $\overline{\mathbb C}$. Hence
the set $\overline{\mathbb C}\setminus C\ne\emptyset$ is also open
in $\overline{\mathbb C}$. The latter contradicts the connectivity
of $\overline{\mathbb C}$, see e.g. Proposition I.1.1 in \cite{FL}.

The obtained contradiction disproves the assumption that $c\ne 0$.
Thus, we have proved that $f$ is extended to a homeomorphism of
$\overline{\mathbb C}$ onto itself with $f(\infty)=\infty$.
\end{proof} $\Box$

\bigskip

Choosing $\psi_{z_0,\varepsilon}(t)\equiv 1/\left(t\,
\log\left(1/t\right)\right)$ in Lemma 3, we obtain by Proposition 2
the following.

\medskip

{\bf Theorem 5.} {\it Let $\mu : \mathbb C\to{\Bbb C}$ be measurable
with $|\mu (z)| < 1$  a.e., $K_{\mu}\in L^1_{\rm loc}(\mathbb C)$
and
\begin{equation}\label{eqINFTYreRE}
\int\limits_{R_0<|z|<R}\ K_{{\mu}}(z)\ {\psi}^2(|z|)\
\frac{dm(z)}{|z|^4}\ =\ o(I^2(R)) \ \ \ \ \ \ \hbox{ as $R\to\infty
$}
\end{equation}
for some $R_0>0$ and a measurable function
${\psi}:(0,\infty)\to(0,\infty)$ such that
\begin{equation}\label{eqIreRE}
I(R)\ \colon =\ \int\limits_{R_0}^{R}{\psi}(t)\ \frac{dt}{t^2}\ <\
\infty\ \ \ \ \ \ \forall\ R\in(R_0,\infty)\ .\end{equation} Suppose
also that $K^T_{\mu}(z,z_0)\leqslant Q_{z_0}(z)$ a.e. in $U_{z_0}$
for every point $z_0\in \mathbb C$, a neighborhood $U_{z_0}$ of
$z_0$ and a function $Q_{z_0}: U_{z_0}\to[0,\infty]$ in the class
${\rm FMO}({z_0})$.

Then the Beltrami equation (\ref{eqBeltrami}) has a regular
homeomorphic solution $f$ in $\mathbb C$ with the
nor\-ma\-li\-za\-tion $f(0)=0$, $f(1)=1$ and $f(\infty)=\infty$.}

\medskip

In particular, by Proposition 1 the conclusion of Theorem 5 holds if
every point $z_0\in\mathbb C$ is the Lebesgue point of the function
$Q_{z_0}$.

\bigskip

By Corollary 1 we obtain the next nice consequence of Theorem 5,
too.

\medskip

{\bf Corollary 8.} {\it Let $\mu : \mathbb C\to{\Bbb C}$ be
measurable with $|\mu (z)| < 1$  a.e., $K_{\mu}\in L^1_{\rm
loc}(\mathbb C)$, (\ref{eqINFTYreRE}) and
\begin{equation}\label{eqMEAN}\overline{\lim\limits_{\varepsilon\to0}}\quad
\dashint_{\mathbb
D(z_0,\varepsilon)}K^T_{\mu}(z,z_0)\,dm(z)<\infty\qquad\forall\
z_0\in\mathbb C\, .\end{equation} Then the Beltrami equation
(\ref{eqBeltrami}) has a regular homeomorphic solution $f$ in
$\mathbb C$ with the nor\-ma\-li\-za\-tion $f(0)=0$, $f(1)=1$ and
$f(\infty)=\infty$.}

\medskip

By (\ref{eqConnect}), we also obtain the following consequences of
Theorem 5.

\medskip

{\bf Corollary 9.} {\it Let $\mu : \mathbb C\to{\Bbb C}$ be
measurable with $|\mu (z)| < 1$  a.e., (\ref{eqINFTYreRE}) and
$K_{\mu}$ have a dominant $Q:\mathbb C\to[1,\infty)$ in the class
{\rm BMO}$_{\rm loc}$. Then the Beltrami equation (\ref{eqBeltrami})
has a regular homeomorphic solution $f$ in $\mathbb C$ with the
nor\-ma\-li\-za\-tion $f(0)=0$, $f(1)=1$ and $f(\infty)=\infty$.}

\medskip

{\bf Remark 8.} In particular, the conclusion of Corollary 7 holds
if $Q\in{\rm W}^{1,2}_{\rm loc}$ because $W^{\,1,2}_{\rm loc}
\subset {\rm VMO}_{\rm loc}$, see e.g. \cite{BN}.

\medskip

{\bf Corollary 10.} {\it Let $\mu : \mathbb C\to{\Bbb C}$ be
measurable with $|\mu (z)| < 1$  a.e., (\ref{eqINFTYreRE}) and
$K_{\mu}(z)\leqslant Q(z)$ a.e. in $\mathbb C$ with a function $Q$
in the class ${\rm FMO}(\mathbb C)$. Then the Beltrami equation
(\ref{eqBeltrami}) has a regular homeomorphic solution $f$ in
$\mathbb C$ with the nor\-ma\-li\-za\-tion $f(0)=0$, $f(1)=1$ and
$f(\infty)=\infty$.}

\medskip

Similarly, choosing  $\psi_{z_0,\varepsilon}(t)\equiv 1/t$ in Lemma
3, we come to the next statement.

\medskip

{\bf Theorem 6.} {\it Let $\mu : \mathbb C\to{\Bbb C}$ be measurable
with $|\mu (z)| < 1$  a.e., $K_{\mu}\in L^1_{\rm loc}(\mathbb C)$,
(\ref{eqINFTYreRE}) and
\begin{equation}\label{eqLOG}
\int\limits_{\varepsilon<|z-z_0|<\varepsilon_0}K^T_{\mu}(z,z_0)\,\frac{dm(z)}{|z-z_0|^2}
=o\left(\left[\log\frac{1}{\varepsilon}\right]^2\right)\qquad\hbox{as
$\varepsilon\to 0$}\qquad\forall\ z_0\in\mathbb C\end{equation} for
some $\varepsilon_0=\varepsilon(z_0)>0$. Then the Beltrami equation
(\ref{eqBeltrami}) has a regular homeomorphic solution $f$ in
$\mathbb C$ with the nor\-ma\-li\-za\-tion $f(0)=0$, $f(1)=1$ and
$f(\infty)=\infty$.}

\medskip

{\bf Remark 9.} Choosing $\psi_{z_0,\varepsilon}(t)\equiv
1/(t\log{1/t})$ instead of $\psi(t)=1/t$ in Lemma 2, we are able to
replace (\ref{eqLOG}) by
\begin{equation}\label{eqLOGLOG}
\int\limits_{\varepsilon<|z-z_0|<\varepsilon_0}\frac{K^T_{\mu}(z,z_0)\,dm(z)}
{\left(|z-z_0|\log{\frac{1}{|z-z_0|}}\right)^2}
=o\left(\left[\log\log\frac{1}{\varepsilon}\right]^2\right)\end{equation}
In general, we are able to give here the whole scale of the
corresponding conditions in $\log$ using functions $\psi(t)$ of the
form
$1/(t\log{1}/{t}\cdot\log\log{1}/{t}\cdot\ldots\cdot\log\ldots\log{1}/{t})$.

\medskip

Now, choosing in Lemma 3 the functional parameter
${\psi}_{z_0,{\varepsilon}}(t) \equiv {\psi}_{z_0}(t) \colon =
1/[tk^T_{\mu}(z_0,t)]$, where $k_{\mu}^T(z_0,r)$ is the integral
mean value of $K^T_{{\mu}}(z,z_0)$ over the circle $S(z_0,r)\, :=\,
\{ z \in\mathbb C:\, |z-z_0|\,=\, r\}$, we obtain one more important
conclusion.

\medskip

{\bf Theorem 7.} {\it Let $\mu : \mathbb C\to{\Bbb C}$ be measurable
with $|\mu (z)| < 1$  a.e., $K_{\mu}\in L^1_{\rm loc}(\mathbb C)$,
(\ref{eqINFTYreRE}) and
\begin{equation}\label{eqLEHTO}\int\limits_{0}^{\varepsilon_0}
\frac{dr}{rk^T_{\mu}(z_0,r)}=\infty\qquad\forall\ z_0\in\mathbb
C\end{equation} for some $\varepsilon_0=\varepsilon(z_0)>0$. Then
the Beltrami equation (\ref{eqBeltrami}) has a regular homeomorphic
solution $f$ in $\mathbb C$ with the nor\-ma\-li\-za\-tion $f(0)=0$,
$f(1)=1$ and $f(\infty)=\infty$.}

\medskip

{\bf Corollary 11.} {\it Let $\mu : \mathbb C\to{\Bbb C}$ be
measurable with $|\mu (z)| < 1$  a.e., $K_{\mu}\in L^1_{\rm
loc}(\mathbb C)$, (\ref{eqINFTYreRE}) and
\begin{equation}\label{eqLOGk}k^T_{\mu}(z_0,\varepsilon)=O\left(\log\frac{1}{\varepsilon}\right)
\qquad\mbox{as}\ \varepsilon\to0\qquad\forall\ z_0\in\mathbb C\
.\end{equation} Then the Beltrami equation (\ref{eqBeltrami}) has a
regular homeomorphic solution $f$ in $\mathbb C$ with the
nor\-ma\-li\-za\-tion $f(0)=0$, $f(1)=1$ and $f(\infty)=\infty$. }

\medskip

{\bf Remark 10.} In particular, the conclusion of Corollary 10 holds
if
\begin{equation}\label{eqLOGK} K^T_{\mu}(z,z_0)=O\left(\log\frac{1}{|z-z_0|}\right)\qquad{\rm
as}\quad z\to z_0\quad\forall\ z_0\in\overline{D}\,.\end{equation}
Moreover, the condition (\ref{eqLOGk}) can be replaced by the whole
series of more weak conditions
\begin{equation}\label{edLOGLOGk}
k^T_{\mu}(z_0,\varepsilon)=O\left(\left[\log\frac{1}{\varepsilon}\cdot\log\log\frac{1}
{\varepsilon}\cdot\ldots\cdot\log\ldots\log\frac{1}{\varepsilon}
\right]\right) \qquad\forall\ z_0\in \overline{D}\ .
\end{equation}

\medskip

For further consequences, the following statement is useful, see
e.g. Theorem 3.2 in \cite{RSY$_5$}.

\medskip

{\bf Proposition 7.} {\it Let $Q:{\Bbb D}\to[0,\infty]$ be a
measurable function such that
\begin{equation}\label{eq5.555} \int\limits_{\Bbb
D}\Phi(Q(z))\,dm(z)<\infty\end{equation} where
$\Phi:[0,\infty]\to[0,\infty]$ is a non-decreasing convex function
such that \begin{equation}\label{eq3.333a}
\int\limits_{\delta}^{\infty}\frac{d\tau}{\tau\Phi^{-1}(\tau)}=\infty\end{equation}
for some $\delta>\Phi(+0)$. Then \begin{equation}\label{eq3.333A}
\int\limits_0^1\frac{dr}{rq(r)}=\infty\end{equation} where $q(r)$ is
the average of the function $Q(z)$ over the circle $|z|=r$.}

\medskip

Here we use the following notions of the inverse function for
monotone functions. Namely, for every non-decreasing function
$\Phi:[0,\infty]\rightarrow[0,\infty]$ the inverse function
$\Phi^{-1}:[0,\infty]\rightarrow[0,\infty]$ can be well-defined by
setting
\begin{equation}\label{eqINVERSE}
\Phi^{-1}(\tau)\ :=\ \inf\limits_{\Phi(t)\geqslant\tau} t\,.
\end{equation}
Here $\inf$ is equal to $\infty$ if the set of $t\in[0,\infty]$ such
that $\Phi(t)\geqslant\tau$ is empty. Note that the function
$\Phi^{-1}$ is non-decreasing, too. It is evident immediately by the
definition that $\Phi^{-1}(\Phi(t)) \leqslant t$ for all
$t\in[0,\infty]$ with the equality except intervals of constancy of
the function $\Phi(t)$.

\medskip

Let us recall the connection of condition (\ref{eq3.333a}) with
other integral conditions, see e.g. Theorem 2.5 in \cite{RSY$_5$}.

\medskip

{\bf Remark 11.} Let $\Phi:[0,\infty]\to[0,\infty]$ be a
non-decreasing function and set
\begin{equation}\label{eqLOGFi}
H(t)\ =\ \log\Phi(t)\ .
\end{equation}
Then the equality
\begin{equation}\label{eq333Frer}\int\limits_{\Delta}^{\infty}H'(t)\,\frac{dt}{t}=\infty,
\end{equation}
implies the equality
\begin{equation}\label{eq333F}\int\limits_{\Delta}^{\infty}
\frac{dH(t)}{t}=\infty\,,\end{equation} and (\ref{eq333F}) is
equivalent to
\begin{equation}\label{eq333B}
\int\limits_{\Delta}^{\infty}H(t)\,\frac{dt}{t^2}=\infty\,\end{equation}
for some $\Delta>0$, and (\ref{eq333B}) is equivalent to each of the
equalities
\begin{equation}\label{eq333C} \int\limits_{0}^{\delta_*}H\left(\frac{1}{t}\right)\,{dt}=\infty\end{equation} for
some $\delta_*>0$, \begin{equation}\label{eq333D}
\int\limits_{\Delta_*}^{\infty}\frac{d\eta}{H^{-1}(\eta)}=\infty\end{equation}
for some $\Delta_*>H(+0)$ and to (\ref{eq3.333a}) for some
$\delta>\Phi(+0)$.

Moreover, (\ref{eq333Frer}) is equivalent to (\ref{eq333F}) and
hence to (\ref{eq333B})–(\ref{eq333D}) as well as to
(\ref{eq3.333a}) are equivalent to each other if $\Phi$ is in
addition absolutely continuous. In particular, all the given
conditions are equivalent if $\Phi$ is convex and non-decreasing.


Note that the integral in (\ref{eq333F}) is understood as the
Lebesgue--Stieltjes integral and the integrals in (\ref{eq333Frer})
and (\ref{eq333B})--(\ref{eq333D}) as the ordinary Lebesgue
integrals. It is necessary to give one more explanation. From the
right hand sides in the conditions (\ref{eq333Frer})--(\ref{eq333D})
we have in mind $+\infty$. If $\Phi(t)=0$ for $t\in[0,t_*]$, then
$H(t)=-\infty$ for $t\in[0,t_*]$ and we complete the definition
$H'(t)=0$ for $t\in[0,t_*]$. Note, the conditions (\ref{eq333F}) and
(\ref{eq333B}) exclude that $t_*$ belongs to the interval of
integrability because in the contrary case the left hand sides in
(\ref{eq333F}) and (\ref{eq333B}) are either equal to $-\infty$ or
indeterminate. Hence we may assume in
(\ref{eq333Frer})--(\ref{eq333C}) that $\delta >t_0$,
correspondingly, $\Delta<1/t_0$ where
$t_0:=\sup\limits_{\Phi(t)=0}t$, and set $t_0=0$ if $\Phi(0)>0$.


The most interesting of the above conditions is (\ref{eq333B}) that
can be rewritten in the form:
\begin{equation}\label{eq5!}
\int\limits_{\Delta}^{\infty}\log\, \Phi(t)\ \frac{dt}{t^{2}}\ =\
+\infty\ \ \ \ \ \ \mbox{for some $\Delta > 0$}\ .
\end{equation}

\medskip

Combining Theorems 7, Proposition 7 and Remark 11, we obtain the
following result.

\medskip

{\bf Theorem 8.} {\it Let $\mu : \mathbb C\to{\Bbb C}$ be measurable
with $|\mu (z)| < 1$  a.e., $K_{\mu}\in L^1_{\rm loc}(\mathbb C)$,
(\ref{eqINFTYreRE}) and
\begin{equation}\label{eqINTEGRAL}\int\limits_{U_{z_0}}\Phi_{z_0}\left(K^T_{\mu}(z,z_0)\right)\,dm(z)<\infty
\qquad\forall\ z_0\in \mathbb C\end{equation} for a neighborhood
$U_{z_0}$ of $z_0$ and a convex non-decreasing function
$\Phi_{z_0}:[0,\infty]\to[0,\infty]$ with
\begin{equation}\label{eqINT}
\int\limits_{\Delta(z_0)}^{\infty}\log\,\Phi_{z_0}(t)\,\frac{dt}{t^2}\
=\ +\infty\end{equation} for some $\Delta(z_0)>0$. Then the Beltrami
equation (\ref{eqBeltrami}) has a regular homeomorphic solution $f$
in $\mathbb C$ with the nor\-ma\-li\-za\-tion $f(0)=0$, $f(1)=1$ and
$f(\infty)=\infty$.}

\medskip

{\bf Corollary 12.} {\it Let $\mu : \mathbb C\to{\Bbb C}$ be
measurable with $|\mu (z)| < 1$  a.e., $K_{\mu}\in L^1_{\rm
loc}(\mathbb C)$, (\ref{eqINFTYreRE}) and
\begin{equation}\label{eqEXP}\int\limits_{U_{z_0}}e^{\alpha(z_0) K^T_{\mu}(z,z_0)}\,dm(z)<\infty
\qquad\forall\ z_0\in \mathbb C\end{equation} for some
$\alpha(z_0)>0$ and a neighborhood $U_{z_0}$ of the point $z_0$.
Then the Beltrami equation (\ref{eqBeltrami}) has a regular
homeomorphic solution $f$ in $\mathbb C$ with the
nor\-ma\-li\-za\-tion $f(0)=0$, $f(1)=1$ and $f(\infty)=\infty$.}

\medskip

Since $K^T_{\mu}(z,z_0) \leqslant K_{\mu}(z)$ for $z$ and $z_0\in
\Bbb C$, we also obtain the following consequences of Theorem 8.

\medskip

{\bf Corollary 13.} {\it Let $\mu : \mathbb C\to{\Bbb C}$ be
measurable with $|\mu (z)| < 1$  a.e., (\ref{eqINFTYreRE}) and
\begin{equation}\label{eqINTK}\int\limits_{C}\Phi\left(K_{\mu}(z)\right)\,dm(z)<\infty\end{equation}
over each compact $C$ in $\mathbb C$ for a convex non-decreasing
function $\Phi:[0,\infty]\to[0,\infty]$ with
\begin{equation}\label{eqINTF}
\int\limits_{\delta}^{\infty}\log\,\Phi(t)\,\frac{dt}{t^2}\ =\
+\infty\end{equation} for some $\delta>0$. Then the Beltrami
equation (\ref{eqBeltrami}) has a regular homeomorphic solution $f$
in $\mathbb C$ with the nor\-ma\-li\-za\-tion $f(0)=0$, $f(1)=1$ and
$f(\infty)=\infty$.}

\medskip

{\bf Corollary 14.} {\it Let $\mu : \mathbb C\to{\Bbb C}$ be
measurable with $|\mu (z)| < 1$  a.e., (\ref{eqINFTYreRE}) and, for
some $\alpha>0$, over each compact $C$ in $\mathbb C$,
\begin{equation}\label{eqEXPA}\int\limits_{C}e^{\alpha K_{\mu}(z)}\,dm(z)\ <\
\infty\ .
\end{equation} Then the Beltrami
equation (\ref{eqBeltrami}) has a regular homeomorphic solution $f$
in $\mathbb C$ with the nor\-ma\-li\-za\-tion $f(0)=0$, $f(1)=1$ and
$f(\infty)=\infty$.}

\section{On existence of solutions with asymptotics at infinity}

In the extended complex plane $\overline{\mathbb C}=\mathbb C\cup
\{\infty\}$, we will use the so-called {\bf spherical area} whose
element can be given through the element $dm(z)$ of the Lebesgue
measure (usual area)
\begin{equation}\label{eqSAREA}
dS(z)\ :=\ \frac{4\,d\,m(z)}{(1+|z|^2)^2}\ =\
\frac{4\,dx\,dy}{(1+|z|^2)^2}\ ,\ \ \ \ z=x+iy\ .
\end{equation}

Let us start from the following general lemma on the existence of
regular homeomorhic solutions for the Beltrami equations in $\mathbb
C$ with asymptotic homogeneity at infinity.

\bigskip

{\bf Lemma 4.} {\it Let a function $\mu : \mathbb C\to{\Bbb C}$ be
measurable with $|\mu (z)| < 1$  a.e., $K_{\mu}$ have a majorant $Q$
of the class BMO in a connected open (punctured at $\infty$)
neighborhood $U$ of infinity,
\begin{equation}\label{eqILEBEG}
\int\limits_{|z|>R} |\mu(z)|\ dS(z)\ =\ o\left( \frac{1}{R^2}\right)
\end{equation}
and, moreover,
\begin{equation}\label{eqISTRONG}
\int\limits_{|z|>R}K_{\mu}(z)\ dS(z)\ =\ O\left(
\frac{1}{R^2}\right)\ .
\end{equation}

Suppose also that, for every $z_0\in{\mathbb C\setminus U}$, there
exist $\varepsilon_0=\varepsilon(z_0)>0$ and a family of measurable
functions ${\psi}_{z_0,{\varepsilon}}:(0,\infty)\to(0,\infty)$ such
that
\begin{equation}\label{eqIREG}
I_{z_0}({\varepsilon})\ \colon =\
\int\limits_{{\varepsilon}}^{{\varepsilon}_0}{\psi}_{z_0,{\varepsilon}}(t)\
dt\ <\ \infty\ \ \ \ \ \ \forall\
{\varepsilon}\in(0,{\varepsilon}_0)\end{equation} and
\begin{equation}\label{eqTREG}
\int\limits_{{\varepsilon}<|z-z_0|<{\varepsilon}_0}\
K^T_{{\mu}}(z,z_0)\cdot{\psi}^2_{z_0,{\varepsilon}}(|z-z_0|)\ dm(z)\
=\ o(I^2_{z_0}({\varepsilon})) \ \ \  \ \ \ \hbox{ as
${\varepsilon}\to 0$}\ \ \ \forall\ z_0\in \mathbb C\ .
\end{equation}

Then the Beltrami equation (\ref{eqBeltrami}) has an approximate
homeomorphic solution $f$ in $\mathbb C$ with $f(0)=0$, $f(1)=1$ and
$f(\infty)=\infty$ that is asymptotically homogeneous at infinity,
$f(\zeta z)\sim\zeta f(z)$ as $z\to\infty$ for all $\zeta\in\mathbb
C$, i.e.,
\begin{equation}\label{eq_4.7INFTY}
  \boxed{\ \lim\limits_{\underset{z\in\mathbb C}{z \to \infty ,}} \frac{f(z \zeta)}{f(z)} = \zeta\ }
  \ \ \ \ \ \ \ \ \ \ \ \ \ \ \forall\ \zeta \in \mathbb{C}
\end{equation} and the limit (\ref{eq_4.7INFTY}) is locally uniform
with respect to the parameter $\zeta$ in $\mathbb C$.}

\bigskip

{\bf Remark 12.} (\ref{eqILEBEG}) and (\ref{eqISTRONG}) can be
replaced by only one (stronger) condition
\begin{equation}\label{eqISTRONGer}
\lim\limits_{r\to\infty}
\frac{R^2}{\pi}\int\limits_{|z|>R}K_{\mu}(z)\ dS(z)\ =\ 1\ .
\end{equation}

Note also that, arguing similarly to the proofs of Theorem 1 and
Corollary 7, we see that the locally uniform property of the
asymptotic homogeneity of $f$ at infinity (\ref{eq_4.7INFTY})
implies its {\bf conformality by Belinskij at infinity}, i.e.,
\begin{equation}\label{eqDifBelCINF}
\boxed{\ f(z)\ =\ A(\rho)\cdot \left[\ z\ +\ o(\rho)\ \right]\
 }\quad\quad\quad \mbox{as $z\to \infty$}\ ,
\end{equation}
where $A(\rho)$ depends only on $\rho = |z|$, $o(\rho) / \rho \to 0$
as $\rho \to \infty$ and, moreover,
\begin{equation}\label{eqDifBelRINF}
\boxed{\ \lim\limits_{\rho \to \infty} \frac{A(t \rho)}{A(\rho)} =
1\ } \quad\quad\quad
    \forall\,t>0\ ,
\end{equation}
its {\bf conformality by Lavrent'iev  at infinity}, i.e.,
\begin{equation}\label{eq_4.8INFTY}
    \boxed{\ \lim\limits_{R \to \infty}\ \frac{\max\limits_{|z| =R} |f(z)|} {\min\limits_{|z| =R} |f(z)|}\ =\
    1\ }\ ,
\end{equation}
the logarithmic property at infinity
\begin{equation}\label{eqLOGARITHM}
\boxed{\lim\limits_{z\to\infty}\ \frac{\ln\,|f(z)|}{\ln\,|z|}\ =\ 1\
}\ ,
\end{equation}
asymptotic preserving angles at infinity, i.e.,
\begin{equation}\label{eq_4.9INFTY}
    \boxed{\ \lim\limits_{z \to \infty}\ \arg\, \left[\frac{f(z \zeta)}{f(z)} \right] = \arg
    \zeta\ }\ \ \ \ \ \forall\ \zeta \in \mathbb{C}^*
\end{equation}
and asymptotic preserving moduli of rings at infinity, i.e.,
\begin{equation}\label{eq_4.10INFTY}
    \boxed{\ \lim\limits_{z \to \infty}\ \frac{|f(z \,\zeta)|}{|f(z)|}\ =\
    |\zeta|\ }\ \ \ \ \ \forall\ \zeta \in \mathbb{C}^*\ .
\end{equation}
The latter two geometric properties characterize asymptotic
homogeneity at infinity and demonstrate that it is very close to the
usual conformality at infinity.

\begin{proof}
The extended complex plane $\overline{\mathbb C}=\mathbb
C\cup\{\infty\}$ is a metric space with a measure with respect to
the spherical (chordal) metric $s$, see (\ref{eq3C}), and the
spherical area $S$, see (\ref{eqSAREA}). This space is regular by
Ahlfors that is evident from the geometric interpretation of
$\overline{\mathbb C}$ as the so-called stereographic projection of
a sphere in $\mathbb R^3$, see details e.g. in Section 13 and
Supplement B in the monograph \cite{MRSY}.

Let us recall only here that, if the function $Q$ belongs to the
class BMO in $U$ with respect to the Euclidean distance and the
usual area in $\mathbb C$, then $Q$ is in BMO with respect to the
spherical distance and the spherical area not only in $U$ but also
in $U\cup\{\infty\}$, see Lemma B.3 and Proposition B.1 in
\cite{MRSY}. Moreover, we have an analog of Proposition 2 in terms
of spherical metric and area, see Lemma 13.2 and Remark 13.3 in
\cite{MRSY}, that in turn can be rewritten in terms of the Euclidean
distance and area at infinity in the following form:
\begin{equation}\label{eqINFTYreS}
\int\limits_{R_0<|z|<R}\ \frac{Q(z)}{\log^2\,|z|}\
\frac{dm(z)}{|z|^2}\ =\ O(\log\log\, R) \ \ \ \ \ \ \hbox{ as
$R\to\infty $}
\end{equation}
for large enough $R_0$ with $\{ z\in \mathbb C: |z|>R_0\}\subseteq
U$. Consequently, we have the condition (\ref{eqINFTYre}) with
$\psi_R(t)\equiv\psi(t):=t^{-1}log\, t$ and by Lemma 3, see also
Remark 7, the Beltrami equation (\ref{eqBeltrami}) has an
approximate solution $f$ in $\mathbb C$ with the
nor\-ma\-li\-za\-tion $f(0)=0$, $f(1)=1$ and $f(\infty)=\infty$.
Recall that $f$ is its regular homeomorphic solution by Proposition
4.

Setting $f^*(\xi):=1/f(1/\xi)$ in $\overline{\mathbb C}$, we see
that $f^*(0)=0$, $f^*(1)=1$, $f^*(\infty)=\infty$ and that $f^*$ is
an approximate solution in $\mathbb C^*=\mathbb C\setminus\{ 0\}$ of
the Beltrami equation with
\begin{equation}\label{eqA1}
\mu^*(\xi)\ :=\
\mu\left(\frac{1}{\xi}\right)\cdot\frac{\xi^2}{\bar\xi^2}\ ,\ \ \ \
\ K_{\mu^*}(\xi)\ =\ K_{\mu}\left(\frac{1}{\xi}\right)\ ,
\end{equation}
because
\begin{equation}\label{eqA2}
f^*_{\bar\xi}(\xi)\ =\ \frac{1}{\bar\xi^2}\cdot\frac{f_{\bar
z}(\frac{1}{\xi})}{f^2(\frac{1}{\xi})}\ ,\ \ \ \ \ \ \
f^*_{\xi}(\xi)\ =\
\frac{1}{\xi^2}\cdot\frac{f_{z}(\frac{1}{\xi})}{f^2(\frac{1}{\xi})}\
\ \ \ \ \mbox{a.e. in $\mathbb C$} \ ,
\end{equation}
see e.g. Section I.C and the proof of Theorem 3 of Section V.B in
\cite{Ahl}.

Note that $f^*$ belongs to the class $W^{1,1}_{\rm loc}(\mathbb
C^*)$ and, consequently, $f^*$ is ACL (absolutely continuous on
lines) in $\mathbb C$, see e.g. Theorems 1 and 2 of Section 1.1.3
and Theorem of Section 1.1.7 in \cite{Ma}. However, it is not clear
directly from (\ref{eqA2}) whether the derivatives $f^*_{\bar\xi}$
and $f^*_{\xi}$ are integrable in a neighborhood of the origin,
because of the first factors in (\ref{eqA2}). Thus, to prove that
$f^*$ is a regular homeomorphic solution of the Beltrami equation in
$\mathbb C$, it remains to establish the latter fact in another way.

Namely, after the replacements of variables $z \longmapsto \xi :=
1/z$ and $R \longmapsto r := 1/R$, in view of (\ref{eqA1}), the
condition (\ref{eqISTRONG}) can be rewritten in the form
\begin{equation}\label{eqISTRONGre}
\limsup\limits_{r\to 0}\
\frac{1}{r^2}\int\limits_{|\xi|<r}K_{\mu^*}(\xi)\ dm(\xi)\ <\
\infty\ ,
\end{equation}
and the latter implies, in particular, that, for some $r_0\in(0,1]$,
\begin{equation}\label{eqISTRONGreZERO}
\frac{1}{r_0^2}\int\limits_{|\xi|<r_0}K_{\mu^*}(\xi)\ dm(\xi)\ <\
\infty\ ,
\end{equation}
i.e., the dilatation quotient $K_{\mu^*}$ of the given Beltrami
equation is integrable in the disk $\mathbb D(r_0)$.

Now, since $f^*$ is a regular homeomorphism in $\mathbb C^*$, in
particular, its Jacobian $J(\xi)=|f^*_{\xi}|^2 -
|f^*_{\bar\xi}|^2\ne 0$ a.e. and hence $|f^*_{\xi}| -
|f^*_{\bar\xi}|\ne 0$ a.e. as well as $f^*_{\xi}\ne 0$ a.e., the
following identities are also correct a.e.
\begin{equation}\label{eqID}
|f^*_{\xi}(\xi)|\ +\ |f^*_{\bar\xi}(\xi)|\ =\ \left[
\frac{|f^*_{\xi}(\xi)|\ +\ |f^*_{\bar\xi}(\xi)|}{|f^*_{\xi}(\xi)|\
-\ |f^*_{\bar\xi}(\xi)|}\right]^{\frac{1}{2}}\cdot\
J^{\frac{1}{2}}(\xi)\ =\ K^{\frac{1}{2}}_{\mu^*}(\xi)\cdot
J^{\frac{1}{2}}(\xi)\ .
\end{equation}
Hence by the H\"older inequality for integrals, see e.g. Theorem 189
in \cite{HLP}, we have that
\begin{equation}\label{eqIDEST}
\int\limits_{|\xi|<r_0} \left(|f^*_{\xi}(\xi)|\ +\
|f^*_{\bar\xi}(\xi)|\right)\, dm(\xi)\ \leq\
\left(\int\limits_{|\xi|<r_0}K_{\mu^*}(\xi)\,
dm(\xi)\right)^{\frac{1}{2}}\cdot
\left(\int\limits_{|\xi|<r_0}J(\xi)\, dm(\xi)\right)^{\frac{1}{2}}
\end{equation}
and, since the latter factor in (\ref{eqIDEST}) is estimated by the
area of $f^*(\mathbb D(r_0))$, see e.g. the Lebesgue theorem in
Section III.2.3 of the monograph \cite{LV}, we conclude that both
partial derivatives $f^*_{\xi}$ and $f^*_{\bar\xi}$  are integrable
in the disk $\mathbb D(r_0)$.

Next, note that the function $Q_*(\xi):=Q(1/\xi)$ is of the class
BMO in a neighborhood of the origin  with respect to the spherical
area as well as with respect to the usual area, see e.g. again Lemma
B.3 in \cite{MRSY}, because also the spherical area is invariant
under rotations of the sphere $\mathbb S^2$ in the stereographic
projection. Moreover, by (\ref{eqILEBEG}) and (\ref{eqA1}), we
obtain that
\begin{equation}\label{eqILEBEGC}
\lim\limits_{r\to 0}\ \frac{1}{r^2}\int\limits_{|\xi|<r}
|\mu^*(\xi)|\ dm(\xi)\ =\ 0\ .
\end{equation}
Thus, by Theorems 3 we conclude that $f^*$ is asymptotically
homogeneous at the origin, i.e.,
\begin{equation}\label{eq_5.7AS}
   \lim\limits_{\underset{\xi \in \mathbb{C}^{*}}{\xi \to 0,}} \frac{f^*(\xi \zeta)}{f^*(\xi)} =
   \zeta\ \ \ \ \ \ \forall\ \zeta\in\mathbb C
\end{equation}
and, furthermore, the limit in(\ref{eq_5.7AS}) is locally uniform in
the parameter $\zeta$.

After the inverse replacements of the variables $\xi\longmapsto
w:=1/\xi$ and the functions $f^*(\xi)\longmapsto f(w)=1/f^*(1/w)$
the relation (\ref{eq_5.7AS}) can be rewritten in the form
\begin{equation}\label{eq_5.7ASINV}
   \lim\limits_{\underset{w \in \mathbb{C}}{w \to \infty ,}} \frac{f(w)}{f\left({w}{\zeta}^{-1}\right)} =
   \zeta\ \ \ \ \ \ \forall\ \zeta\in\mathbb C\ .
\end{equation}
Finally, after one more change of variables $w\longmapsto
z:=w{\zeta}^{-1}$, the latter is transformed into
(\ref{eq_4.7INFTY}), where the limit is locally uniform with respect
to the parameter $\zeta\in\mathbb C$.
\end{proof} $\Box$

\bigskip

Choosing $\psi_{z_0,\varepsilon}(t)\equiv 1/\left(t\,
\log\left(1/t\right)\right)$ in Lemma 4, we obtain by Proposition 2
the following.

\medskip

{\bf Theorem 9.} {\it Let a function $\mu : \mathbb C\to{\Bbb C}$ be
measurable with $|\mu (z)| < 1$  a.e., $K_{\mu}$ have a majorant $Q$
of the class BMO in a neighborhood $U$ of $\infty$ and satisfy
(\ref{eqISTRONGer}). Suppose also that $K^T_{\mu}(z,z_0)\leqslant
Q_{z_0}(z)$ a.e. in $U_{z_0}$ for every point $z_0\in \mathbb
C\setminus U$, a neighborhood $U_{z_0}$ of $z_0$ and a function
$Q_{z_0}: U_{z_0}\to[0,\infty]$ in the class ${\rm FMO}({z_0})$.
Then the Beltrami equation (\ref{eqBeltrami}) has a regular
homeomorphic solution $f$ in $\mathbb C$ with $f(0)=0$, $f(1)=1$ and
$f(\infty)=\infty$ that is asymptotically homogeneous at infinity.}

\medskip

As a particular case of Theorem 9, we obtain the following central
theorem in terms of BMO.

\medskip

{\bf Theorem 10.} {\it Let a function $\mu : \mathbb C\to{\Bbb C}$
be measurable with $|\mu (z)| < 1$  a.e., $K_{\mu}$ have a majorant
$Q$ of the class {\rm BMO(}$\mathbb C${\rm )} and satisfy
(\ref{eqISTRONGer}). Then the Beltrami equation (\ref{eqBeltrami})
has a regular homeomorphic solution $f$ in $\mathbb C$ with
$f(0)=0$, $f(1)=1$ and $f(\infty)=\infty$ that is asymptotically
homogeneous at infinity.}

\medskip

Note also that, in particular, by Proposition 1 the conclusion of
Theorem 9 holds if every point $z_0\in\mathbb C\setminus U$ is the
Lebesgue point of the function $Q_{z_0}$.

\bigskip

By Corollary 1 we obtain the next fine consequence of Theorem 9,
too.

\medskip

{\bf Corollary 15.} {\it Let $\mu : \mathbb C\to{\Bbb C}$ be a
measurable function with $|\mu (z)| < 1$  a.e., $K_{\mu}$ have a
majorant $Q$ of the class BMO in a neighborhood $U$ of $\infty$,
satisfy (\ref{eqISTRONGer}) and
\begin{equation}\label{eqMEANAS}\overline{\lim\limits_{\varepsilon\to0}}\quad
\dashint_{\mathbb
D(z_0,\varepsilon)}K^T_{\mu}(z,z_0)\,dm(z)<\infty\qquad\forall\
z_0\in\mathbb C\setminus U\, .\end{equation} Then the Beltrami
equation (\ref{eqBeltrami}) has a regular homeomorphic solution $f$
in $\mathbb C$ with $f(0)=0$, $f(1)=1$ and $f(\infty)=\infty$ that
is asymptotically homogeneous at infinity.}

\medskip

By (\ref{eqConnect}), we also obtain the following consequences of
Theorem 9.

\medskip

{\bf Corollary 16.} {\it Let $\mu : \mathbb C\to{\Bbb C}$ be a
measurable function with $|\mu (z)| < 1$  a.e.,  $K_{\mu}$ have a
majorant $Q$ of the class BMO in a neighborhood $U$ of $\infty$,
satisfy (\ref{eqISTRONGer}) and $K_{\mu}$ have a dominant
$Q_*:\mathbb C\setminus U\to[1,\infty)$ in the class {\rm BMO}$_{\rm
loc}$. Then the Beltrami equation (\ref{eqBeltrami}) has a regular
homeomorphic solution $f$ in $\mathbb C$ with $f(0)=0$, $f(1)=1$ and
$f(\infty)=\infty$ that is asymptotically homogeneous at infinity.}

\medskip

{\bf Remark 13.} In particular, the conclusion of Corollary 14 holds
if $Q_*\in{\rm W}^{1,2}_{\rm loc}$ because $W^{\,1,2}_{\rm loc}
\subset {\rm VMO}_{\rm loc}$, see e.g. \cite{BN}.

\medskip

{\bf Corollary 17.} {\it Let $\mu : \mathbb C\to{\Bbb C}$ be a
measurable function with $|\mu (z)| < 1$  a.e.,  $K_{\mu}$ have a
majorant $Q$ of the class BMO in a neighborhood $U$ of $\infty$,
satisfy (\ref{eqISTRONGer}) and $K_{\mu}(z)\leqslant Q_*(z)$ a.e. in
$\mathbb C\setminus U$ with a function $Q:\mathbb C\to\mathbb R^+$
of the class ${\rm FMO}(\mathbb C\setminus U)$. Then the Beltrami
equation (\ref{eqBeltrami}) has a regular homeomorphic solution $f$
in $\mathbb C$ with $f(0)=0$, $f(1)=1$ and $f(\infty)=\infty$ that
is asymptotically homogeneous at infinity.}

\medskip

Similarly, choosing  $\psi_{z_0,\varepsilon}(t)\equiv 1/t$ in Lemma
4, we come also to the next statement.

\medskip

{\bf Theorem 11.} {\it Let $\mu : \mathbb C\to{\Bbb C}$ be a
measurable function with $|\mu (z)| < 1$  a.e., $K_{\mu}$ have a
majorant $Q$ of the class BMO in a neighborhood $U$ of $\infty$,
satisfy (\ref{eqISTRONGer}) and,  for some
$\varepsilon_0=\varepsilon(z_0)>0$,
\begin{equation}\label{eqLOGAS}
\int\limits_{\varepsilon<|z-z_0|<\varepsilon_0}K^T_{\mu}(z,z_0)\,\frac{dm(z)}{|z-z_0|^2}
=o\left(\left[\log\frac{1}{\varepsilon}\right]^2\right)\qquad\hbox{as
$\varepsilon\to 0$}\qquad\forall\ z_0\in\mathbb C\setminus U\
.\end{equation} Then the Beltrami equation (\ref{eqBeltrami}) has a
regular homeomorphic solution $f$ in $\mathbb C$ with $f(0)=0$,
$f(1)=1$ and $f(\infty)=\infty$ that is asymptotically homogeneous
at infinity.}

\medskip

{\bf Remark 14.} Choosing $\psi_{z_0,\varepsilon}(t)\equiv
1/(t\log{1/t})$ instead of $\psi(t)=1/t$ in Lemma 4, we are able to
replace (\ref{eqLOGAS}) by
\begin{equation}\label{eqLOGLOGAS}
\int\limits_{\varepsilon<|z-z_0|<\varepsilon_0}\frac{K^T_{\mu}(z,z_0)\,dm(z)}
{\left(|z-z_0|\log{\frac{1}{|z-z_0|}}\right)^2}
=o\left(\left[\log\log\frac{1}{\varepsilon}\right]^2\right)\end{equation}
In general, we are able to give here the whole scale of the
corresponding conditions in $\log$ using functions $\psi(t)$ of the
form
$1/(t\log{1}/{t}\cdot\log\log{1}/{t}\cdot\ldots\cdot\log\ldots\log{1}/{t})$.

\medskip

Now, choosing in Lemma 4 the functional parameter
${\psi}_{z_0,{\varepsilon}}(t) \equiv {\psi}_{z_0}(t) \colon =
1/[tk^T_{\mu}(z_0,t)]$, where $k_{\mu}^T(z_0,r)$ is the average of
$K^T_{{\mu}}(z,z_0)$ over the circle $S(z_0,r)\, :=\, \{ z
\in\mathbb C:\, |z-z_0|\,=\, r\}$, we obtain one more important
conclusion.

\medskip

{\bf Theorem 12.} {\it Let $\mu : \mathbb C\to{\Bbb C}$ be a
measurable function with $|\mu (z)| < 1$  a.e.,  $K_{\mu}$ have a
majorant $Q$ of the class BMO in a neighborhood $U$ of $\infty$,
satisfy (\ref{eqISTRONGer}) and,  for some
$\varepsilon_0=\varepsilon(z_0)>0$,
\begin{equation}\label{eqLEHTOAS}\int\limits_{0}^{\varepsilon_0}
\frac{dr}{rk^T_{\mu}(z_0,r)}=\infty\qquad\forall\ z_0\in\mathbb
C\setminus U\ .\end{equation} Then the Beltrami equation
(\ref{eqBeltrami}) has a regular homeomorphic solution $f$ in
$\mathbb C$ with $f(0)=0$, $f(1)=1$ and $f(\infty)=\infty$ that is
asymptotically homogeneous at infinity.}

\medskip

{\bf Corollary 18.} {\it Let $\mu : \mathbb C\to{\Bbb C}$ be a
measurable function with $|\mu (z)| < 1$  a.e.,  $K_{\mu}$ have a
majorant $Q$ of the class BMO in a neighborhood $U$ of $\infty$,
satisfy (\ref{eqISTRONGer}) and
\begin{equation}\label{eqLOGkAS}k^T_{\mu}(z_0,\varepsilon)=O\left(\log\frac{1}{\varepsilon}\right)
\qquad\mbox{as}\ \varepsilon\to0\qquad\forall\ z_0\in\mathbb
C\setminus U\ .\end{equation} Then the Beltrami equation
(\ref{eqBeltrami}) has a regular homeomorphic solution $f$ in
$\mathbb C$ with $f(0)=0$, $f(1)=1$ and $f(\infty)=\infty$ that is
asymptotically homogeneous at infinity.}

\medskip

{\bf Remark 15.} In particular, the conclusion of Corollary 18 holds
if
\begin{equation}\label{eqLOGKAS} K^T_{\mu}(z,z_0)=O\left(\log\frac{1}{|z-z_0|}\right)\qquad{\rm
as}\quad z\to z_0\quad\forall\ z_0\in\mathbb C\setminus
U\,.\end{equation} Moreover, the condition (\ref{eqLOGkAS}) can be
replaced by the whole series of more weak conditions
\begin{equation}\label{edLOGLOGkAS}
k^T_{\mu}(z_0,\varepsilon)=O\left(\left[\log\frac{1}{\varepsilon}\cdot\log\log\frac{1}
{\varepsilon}\cdot\ldots\cdot\log\ldots\log\frac{1}{\varepsilon}
\right]\right) \qquad\forall\ z_0\in\mathbb C\setminus U\ .
\end{equation}

\medskip

Combining Theorems 12, Proposition 4 and Remark 1, we obtain the
following result.

\medskip

{\bf Theorem 13.} {\it Let $\mu : \mathbb C\to{\Bbb C}$ be a
measurable function with $|\mu (z)| < 1$  a.e., $K_{\mu}$ have a
majorant $Q$ of the class BMO in a neighborhood $U$ of $\infty$,
satisfy (\ref{eqISTRONGer}) and
\begin{equation}\label{eqINTEGRALAS}\int\limits_{U_{z_0}}\Phi_{z_0}\left(K^T_{\mu}(z,z_0)\right)\,dm(z)<\infty
\qquad\forall\ z_0\in\mathbb C\setminus U\end{equation} for a
neighborhood $U_{z_0}$ of $z_0$ and a convex non-decreasing function
$\Phi_{z_0}:[0,\infty]\to[0,\infty]$ with
\begin{equation}\label{eqINTAS}
\int\limits_{\Delta(z_0)}^{\infty}\log\,\Phi_{z_0}(t)\,\frac{dt}{t^2}\
=\ +\infty\ \ \ \ \ \mbox{for some $\Delta(z_0)>0$\ .}\end{equation}
Then the Beltrami equation (\ref{eqBeltrami}) has a regular
homeomorphic solution $f$ in $\mathbb C$ with $f(0)=0$, $f(1)=1$ and
$f(\infty)=\infty$ that is asymptotically homogeneous at infinity.}

\medskip

{\bf Corollary 19.} {\it Let $\mu : \mathbb C\to{\Bbb C}$ be a
measurable function with $|\mu (z)| < 1$  a.e.,   $K_{\mu}$ have a
majorant $Q$ of the class BMO in a neighborhood $U$ of $\infty$,
satisfy (\ref{eqISTRONGer}) and, for some $\alpha(z_0)>0$ and a
neighborhood $U_{z_0}$ of the point $z_0$,
\begin{equation}\label{eqEXPAS}\int\limits_{U_{z_0}}e^{\alpha(z_0) K^T_{\mu}(z,z_0)}\,dm(z)<\infty
\qquad\forall\ z_0\in \mathbb C\setminus U\ .\end{equation}  Then
the Beltrami equation (\ref{eqBeltrami}) has a regular homeomorphic
solution $f$ in $\mathbb C$ with $f(0)=0$, $f(1)=1$ and
$f(\infty)=\infty$ that is asymptotically homogeneous at infinity.}

\medskip

Since $K^T_{\mu}(z,z_0) \leqslant K_{\mu}(z)$ for $z$ and $z_0\in
\Bbb C$, we also obtain the following consequences of Theorem 13.

\medskip

{\bf Corollary 20.} {\it Let $\mu : \mathbb C\to{\Bbb C}$ be a
measurable function with $|\mu (z)| < 1$  a.e., $K_{\mu}$ have a
majorant $Q$ of the class BMO in a neighborhood $U$ of $\infty$,
satisfy (\ref{eqISTRONGer}) and
\begin{equation}\label{eqINTKAS}\int\limits_{\mathbb C\setminus U}\Phi\left(K_{\mu}(z)\right)\,dm(z)<\infty\end{equation}
for a convex non-decreasing function $\Phi:[0,\infty]\to[0,\infty]$
such that, for some $\delta>0$,
\begin{equation}\label{eqINTFAS}
\int\limits_{\delta}^{\infty}\log\,\Phi(t)\,\frac{dt}{t^2}\ =\
+\infty\ .\end{equation} Then the Beltrami equation
(\ref{eqBeltrami}) has a regular homeomorphic solution $f$ in
$\mathbb C$ with $f(0)=0$, $f(1)=1$ and $f(\infty)=\infty$ that is
asymptotically homogeneous at infinity.}

\medskip

{\bf Corollary 21.} {\it Let $\mu : \mathbb C\to{\Bbb C}$ be a
measurable function with $|\mu (z)| < 1$  a.e., $K_{\mu}$ have a
majorant $Q$ of the class BMO in a neighborhood $U$ of $\infty$,
satisfy (\ref{eqISTRONGer}) and, for some $\alpha>0$,
\begin{equation}\label{eqEXPAAS}\int\limits_{\mathbb C\setminus U}e^{\alpha K_{\mu}(z)}\,dm(z)\ <\
\infty\ .
\end{equation} Then the Beltrami equation
(\ref{eqBeltrami}) has a regular homeomorphic solution $f$ in
$\mathbb C$ with $f(0)=0$, $f(1)=1$ and $f(\infty)=\infty$ that is
asymptotically homogeneous at infinity.}

\medskip

{\bf Remark 16.} Recall that by Theorem 5.1 in \cite{RSY$_5$} the
condition (\ref{eqINTFAS}) is not only sufficient but also necessary
for the existence of regular homeomorphic solutions for all Beltrami
equations (\ref{eqBeltrami}) with the integral constraints
(\ref{eqINTKAS}), see also Remark 11.

\medskip

Finally, these results can be applied to the fluid mechanics in
strictly anisotropic and inhomogeneous media because the Beltrami
equation is a complex form of the main equation of hydromechanics,
see e.g. Theorem 16.1.6 in \cite{AIM}, that will be published
elsewhere.


{\bf \noindent Vladimir Gutlyanskii, Vladimir Ryazanov} \\
{\bf 1.} Institute of Applied Mathematics and Mechanics\\
of National Academy of Sciences of Ukraine, \\
Slavyansk 84 100,  UKRAINE\\
{\bf 2.}  Bogdan Khmelnytsky National University of Cherkasy,
\\
Dept. of Phys., Lab. of Math. Phys.,
\\ Cherkasy,   UKRAINE\\
vgutlyanskii@gmail.com , vl.ryazanov1@gmail.com ,\\
Ryazanov@nas.gov.ua

\medskip
{\bf \noindent Evgeny Sevost'yanov} \\
{\bf 1.} Zhytomyr Ivan Franko State University,  \\
Zhytomyr 10 008, UKRAINE \\
{\bf 2.} Institute of Applied Mathematics and Mechanics\\
of National Academy of Sciences of Ukraine, \\
Slavyansk 84 100,  UKRAINE\\
esevostyanov2009@gmail.com

{\bf \noindent Eduard Yakubov} \\
H.I.T. - Holon Institute of Technology,\\
Holon, PO Box 305,  ISRAEL\\
yakubov@hit.ac.il


\begin{thebibliography}{100}

\bibitem{Ahl} {\sc Ahlfors, L.:} Lectures on Quasiconformal Mappings. - Van Nostrand, New York, 1966.

\bibitem{And}
{\sc Andreian Cazacu C.:} On the length-area dilatation. - Complex
Variables, Theory Appl. 50:7--11, 2005, 765--776.

\bibitem{AIM} {\sc
Astala K., Iwaniec T., Martin G.:} Elliptic partial differential
equations and quasiconformal mappings in the plane. - Princeton
Mathematical Series 48. Princeton, NJ: Princeton University Press,
2009.

\bibitem{Be}
{\sc Belinskii P.P.:} General Properties of Quasiconformal Mappings
[in Russian], Nauka, Novosibirsk, 1974.

\bibitem{Bojar} {\sc
Bojarski B.:} Generalized solutions of a system of differential
equations of the first order of the elliptic type with discontinuous
coefficients. - Mat. Sb., N. Ser. 43(85):4, 1958, 451–503.

\bibitem{BGR} {\sc Bojarski B., Gutlyanskii V., Ryazanov V.:}
On the Beltrami equations with two characteristics. - Complex Var.
Ellipt. Eq., 54, No. 10, 935–950 (2009).

\bibitem{BGR$_1$} {\sc Bojarski B., Gutlyanskii V., Ryazanov V.:}
On integral conditions for the general Beltrami equations. Complex
Anal. Oper. Theory 5:3, 2011, 835-845.

\bibitem{BGR$_2$} {\sc Bojarski B., Gutlyanskii V., Ryazanov V.:}
On existence and representation of solutions for general degenerate
Beltrami equations. - Complex Var. Elliptic Equ. 59:1, 2014, 67-75.

\bibitem{BN}
{\sc Brezis H., Nirenberg L.:} {Degree theory and BMO. I. Compact
manifolds without boundaries.} Selecta Math. (N.S.) 1:2, 1995,
197--263.

\bibitem{CFL}
{\sc Chiarenza F., Frasca M., Longo P.:} {$W^{2,p}$-solvability of
the Dirichlet problem for nondivergence elliptic equations with VMO
coefficients.} - Trans. Amer. Math. Soc. 336:2, 1993, 841--853.

\bibitem{DS} {\sc Dunford N., Schwartz Jacob T.:} Linear operators. I. General theory. -
Pure and Applied Mathematics. Vol. 7. Interscience Publishers:  New
York and London, 1958.

\bibitem{EG}
{\sc Evans L.C., Gapiery R.F.:} Measure Theory and Fine Properties
of Functions. - CRC Press, Boca Raton, FL, 1992.

\bibitem{FL} {\sc Fischer W., Lieb I.:} A course in complex analysis. From
basic results to advanced topics. - Vieweg + Teubner: Wiesbaden,
2012.

\bibitem{GeLe}
{\sc Gehring F.W., Lehto O.:} On the total differentiability of
functions of a complex variable. - Ann. Acad. Sci. Fenn. A1. Math.
\textbf{272}, 9 (1959).

\bibitem{GMRV} {\sc Gutlyanskii V.Ya., Martio O., Ryazanov
V.I., Vuorinen M.:} On convergence theorems for space quasiregular
mappings. - Forum Math. 10, no. 3, 353-375 (1998).

\bibitem{GMSV} {\sc
Gutlyanskii V., Martio O., Sugawa T., Vuorinen M.:} On the
degenerate Beltrami equation. - Trans. Amer. Math. Soc. 357:3, 2005,
875-900.

\bibitem{GR}
{\sc Gutlyanskii V.Ya., Ryazanov V.I.:} On the local behaviour of
quasi-conformal mappings. Izv. Math. 59, No. 3, 471-498 (1995).


\bibitem{GRL}
{\sc Gutlyanskii V.Ya., Ryazanov V.I., Lomako T.V.:} To the theory
of variational method for Beltrami equations. - Ukr. Mat. Visn., 8,
No. 4, 513–536 (2011); transl. in J. Math. Sci., 182, No. 1, 37–54
(2012).

\bibitem{GRSY} {\sc Gutlyanskii V., Ryazanov V., Srebro U., Yakubov E.:}
The Beltrami Equation: A Geometric Approach. - Developments in
Mathematics 26, Springer: Berlin, 2012.

\bibitem{GRSY$_*$} {\sc Gutlyanskii V., Ryazanov V., Srebro U., Yakubov E.:}
On recent advances in the Beltrami equations. - J. Math. Sci. (USA)
175:4, 2011, 413-449.

\bibitem{HLP}
{\sc  Hardy G.H., Littlewood J.E., Polya G.:} Inequalities.
Cambridge University Press: Cambridge (UK) etc., 1988.

\bibitem{HKM}
{\sc  Heinonen J., Kilpelainen T., Martio O.:} {Nonlinear Potential
Theory of Degenerate Elliptic Equations.} - Oxford Mathematical
Monographs, Clarendon Press: Oxford - New York - Tokyo, 1993.

\bibitem{IR}
{\sc Ignat'ev A.A., Ryazanov V.I.:} {Finite mean oscillation in the
mapping theory.} - Ukrainian Math. Bull. 2:3, 2005, 403-424.

\bibitem{IS}
{\sc Iwaniec T., C. Sbordone C.:} {Riesz transforms and elliptic
PDEs with VMO coefficients.} - J. Anal. Math. 74, 1998, 183--212.

\bibitem{JN}
{\sc John F., Nirenberg L.:} {On functions of bounded mean
oscillation.} - Comm. Pure Appl. Math. 14, 1961, 415--426.

\bibitem{KR}
{\sc Kolomoitsev Iu.S., Ryazanov V.I.:} Uniqueness of approximate
solutions of Beltrami equations. Tr. Inst. Prikl. Mat. Mekh. 19,
116-124 (2009).

\bibitem{KPRS}
{\sc Kovtonyuk D.A., Petkov I.V., Ryazanov V.I., Salimov R.R.:}
Boundary behavior and the Dirichlet problem for Beltrami equations.
St. Petersbg. Math. J. 25, no. 4, 587-603 (2014); transl. from
Algebra Anal. 25, no. 4, 101-124 (2013).

\bibitem{Ku} {\sc Kuratowski K.:} Topology. Vol. I. New edition,
revised and augmented.(English) Academic Press: New York-London;
PWN-Polish Scientific Publishers: Warszawa, 1966.

\bibitem{LS} {\sc Lawrentjew M.A., Schabat B.W.:} Methoden der
komplexen Funktionentheorie. (German) Mathematik f\"ur
Naturwissenschaft und Technik 13, VEB Deutscher Verlag der
Wissenschaften: Berlin, 1967.

\bibitem{Le}
{\sc Lehto O.:} Homeomorphisms with a prescribed dilatation. -
Lecture Notes in Math. 118, 1968, 58-73.

\bibitem{LV} {\sc Lehto O., Virtanen K.I.:} Quasiconformal mappings in the plane. -
Die Grundlehren der mathematischen Wissenschaften 126, Springer:
Berlin-Heidelberg-New York, 1973.

\bibitem{LSS} {\sc Lomako T., Salimov R., Sevostyanov E.:} On
equicontinuity of solutions to the Beltrami equations. Ann. Univ.
Buchar., Math. Ser. 1(59), no. 2, 263-274 (2010).

\bibitem{MRSY} {\sc Martio O., Ryazanov V., Srebro U., Yakubov E.:} Moduli in modern mapping
theory. - Springer Monographs in Mathematics. - Springer: New York,
2009.

\bibitem{MRV}
{\sc Martio O., Ryazanov V., Vuorinen M.:} {BMO and Injectivity of
Space Quasiregular Mappings.} - Math. Nachr. 205, 1999, 149--161.

\bibitem{Ma} {\sc Maz’ya V.G.:} Sobolev spaces. Berlin
etc.: Springer-Verlag, 1985.

\bibitem{Me}
{\sc Menchoff, D.:} Sur les differentielles totales des fonctions
univalentes. Math. Ann. \textbf{105}, 75--85 (1931).

\bibitem{Pal}
{\sc Palagachev D.K.:} {Quasilinear elliptic equations with VMO
coefficients.} - Trans. Amer. Math. Soc. 347:7, 1995, 2481--2493.

\bibitem{PS}
{\sc Polya G., Szego G.:} Aufgaben und Lehrs\"atze aus der Analysis.
Bd. I: Reihen, Integralrechnung, Funktionentheorie. (German) - Die
Grundlehren der mathematischen Wissenschaften in Einzeldarstellungen
19. Springer: Berlin (1925).

\bibitem{Po}
{\sc Ponomarev S.P.:} The $N^{-1}-$property of mappings, and Lusin’s
$(N)-$condition. Mat. Zametki. – 1995. – 58. – P.411-418; transl. in
Math. Notes. – 1995. – 58. – P. 960-965.

\bibitem{Ra$_1$}
{\sc Ragusa M.A.:} {Elliptic boundary value problem in vanishing
mean oscillation hypothesis.} - Comment. Math. Univ. Carolin. 40:4,
1999, 651--663.

\bibitem{Ra$_2$}
{\sc Ragusa M.A., Tachikawa A.:} Partial regularity of the
minimizers of quadratic functionals with VMO coefficients. J. Lond.
Math. Soc., II. Ser. 72:3, 2005, 609-620.

\bibitem{Rei}
{\sc Reimann H.M.:} Functions of bounded mean oscillation and
quasiconformal mappings. (English) Comment. Math. Helv. 49, 260-276
(1974).

\bibitem{RR}
{\sc Reimann H.M., Rychener T.:} {Funktionen Beschr\"ankter
Mittlerer Oscillation.} - Lecture Notes in Math. 487, 1975.

\bibitem{Re}
{\sc Reshetnyak Yu.G.:} {Space mappings with bounded distortion}. -
Translations of Mathematical Monographs {73}, Providence, RI,
American Mathematical Society (AMS), 1989.

\bibitem{Ru}
{\sc Rudin W.:} Function Theory in Polydiscs. Math. Lect. Notes
Ser., Benjamin, Inc., New-York–Amsterdam (1969).

\bibitem{R}
{\sc Ryazanov V.I.:} A criterion for differentiability in the sense
of Belinskii and its consequences. - Ukr. Math. J. 44, no. 2,
254-258 (1992).

\bibitem{RS}
{\sc Ryazanov R., Salimov R.:} {Weakly flat spaces and boundaries in
the mapping theory.} - Ukrain. Math. Bull. 4:2, 2007, 199--233.

\bibitem{RSS} {\sc Ryazanov V., Salimov R., Sevost'yanov~E.:}
{On Convergence Analysis of Space Homeomorphisms.} - Siberian
Advances in Mathematics. {23,} no. 4, 263--293 (2013).

\bibitem{RSY} {\sc Ryazanov V., Srebro U., Yakubov E.:} BMO-quasiconformal mappings.
- J. d'Anal. Math. 83, 2001, 1-20.

\bibitem{RSY$_1$} {\sc Ryazanov V., Srebro U., Yakubov E.:}
On ring solutions of Beltrami equations. - J. d'Anal. Math. 96,
2005, 117--150.

\bibitem{RSY$_2$}
{\sc Ryazanov V., Srebro U., Yakubov E.:} {Beltrami equation and FMO
functions.} - Contemp. Math. 382, Israel Math. Conf. Proc., 2005,
357--364.

\bibitem{RSY$_3$} {\sc Ryazanov~V., U. Srebro, and E. Yakubov:}
Finite mean oscillation and the Beltrami equation. - Israel Math. J.
153, 2006, 247--266.

\bibitem{RSY$_4$} {\sc Ryazanov~V., U. Srebro, and E. Yakubov:}
On the theory of the Beltrami equation. - Ukr. Math. J. 58:11, 2006,
1786-1798.

\bibitem{RSY$_5$} {\sc Ryazanov V., Srebro U., Yakubov E.:}
Integral conditions in the theory of the Beltrami equations. Complex
Var. Elliptic Equ. 57:12, 2012, 1247-1270.

\bibitem{RSY$_6$}
{\sc Ryazanov V., Srebro U., Yakubov E.:} The  Beltrami equation and
ring homeomorphisms. Ukrainian Math. Bull. {4}, no. 1, 79--115
(2007).

\bibitem{RSY$_8$}  {\sc Ryazanov V., Srebro U., Yakubov E.:}
{On convergence theory for Beltrami equations}. - Укр. мат. вiсник
{5}, no. 4, 524--535 (2008); transl. in Ukr. Math. Bull. {5}, no. 4,
517--528 (2008).

\bibitem{Sal}
{\sc Salimov R.:} On regular homeomorphisms in the plane. Ann. Acad.
Sci. Fenn., Math. 35, no. 1, 285-289 (2010).

\bibitem{Sar}
{\sc Sarason D.:} {Functions of vanishing mean oscillation.} -
Trans. Amer. Math. Soc. 207, 1975, 391--405.

\bibitem{Sp}
{\sc Spanier E.H.:} Algebraic topology. Springer-Verlag: Berlin,
1995.

\bibitem{SY}
{\sc Srebro U. and Yakubov E.:} Beltrami equation - in Handbook of
complex analysis: geometric function theory (ed. K\"uhnau R.), V. 2,
555-597, Elsevier/North Holland, Amsterdam, 2005.

\bibitem{Va} {\sc V\"{a}is\"{a}l\"{a} J.:} {Lectures on $n$-Dimensional
Quasiconformal Mappings}, Lecture Notes in Math. {229}, Berlin etc.,
Springer--Verlag, 1971.


\bibitem{Vek}
{\sc Vekua I.N.:} Generalized analytic functions. - Pergamon Press:
London, 1962.

\end{thebibliography}
\end{document}